\documentclass[preprint,12pt]{elsarticle}

\usepackage{amssymb}
\usepackage{amsthm}
\usepackage{mathrsfs}
\usepackage{enumerate}
\usepackage{amsmath}
\usepackage{amsfonts}
\usepackage{amssymb}
\usepackage{multicol}
\usepackage{bm}
\usepackage{graphicx}
\usepackage{titlesec}
\usepackage{etoolbox}
\usepackage[utf8]{inputenc}
\usepackage[english]{babel}
\usepackage{xcolor,soul}
\sethlcolor{yellow}
\usepackage[framemethod=tikz]{mdframed}
\usepackage{float}
\restylefloat{table}
\usepackage{multirow}
\makeatletter
\patchcmd{\ttlh@hang}{\parindent\z@}{\parindent\z@\leavevmode}{}{}
\patchcmd{\ttlh@hang}{\noindent}{}{}{}
\makeatother
\usepackage{titletoc}
\usepackage{fancyhdr}
\usepackage{indentfirst}
\usepackage[margin=0.8in,vmargin=2.7cm]{geometry} 

\numberwithin{equation}{section}

\journal{arXiv}

\begin{document}
\begin{frontmatter}

\title{Efficient second-order semi-implicit finite element method for fourth-order nonlinear diffusion equations}

\author[addr1]{Sana Keita}     \ead{sana.keita@um6p.ma}
\cortext[corres]{Corresponding author.}
\author[addr1,addr2]{Abdelaziz Beljadid\corref{corres}} \ead{abdelaziz.beljadid@um6p.ma}
\author[addr2]{Yves Bourgault}     \ead{ybourg@uottawa.ca}
\address[addr1]{International Water Research Institute, Mohammed VI Polytechnic University,  Morocco} 
\address[addr2]{Department of Mathematics and Statistics, University of Ottawa, Canada \qquad\qquad\qquad\qquad}

\begin{abstract}
We focus here on a class of fourth-order parabolic equations that can be written as a system of second-order equations by introducing an auxiliary variable. We design a novel second-order fully discrete mixed finite element method to approximate these equations. In our approach, we propose new techniques using the second-order backward differentiation formula for the time derivative and a special technique for the approximation of nonlinear terms. The use of the proposed technique for nonlinear terms makes the developed numerical scheme efficient in terms of computational cost since the proposed method only deals with a linear system at each time step and no iterative resolution is needed. A numerical convergence study is performed using the method of manufactured and analytical solutions of the system where we investigate different boundary conditions. With respect to the spatial discretization, convergence rates are found to at least match a priori error estimates available for linear problems. The convergence analysis is completed with an investigation of the temporal discretization where we numerically demonstrate the second-order time-accuracy of the proposed scheme using the method of reference solution. We present a series of numerical tests to demonstrate the efficiency and robustness of the proposed scheme. 
\end{abstract}


\end{frontmatter}

\section{Introduction}
High order partial differential equations (PDEs) containing the bi-harmonic operator often arise in two-phase and multiphase problems from fluid mechanics and solid mechanics. Examples in the literature include the thin film equations \cite{ehrhard_davis_1991,greenspan1978,Hocking1981,Huppert1982,Moriarty1991}, the Cahn-Hilliard (CH) equations \cite{CAHN1962,CahnHilliard1958}, ocean and atmosphere circulation models \cite{DELHEZ2007}, tumor-growth models \cite{Colli2017,Wise2008}, cancer modeling \cite{Lowengrub_2010}, biological entities \cite{Elson2010,Khain2008}, image processing \cite{Bertozzi4}, the fourth-order model for the approximation of crack propagation \cite{AMIRI2016},  models of unstable infiltration through porous media \cite{Beljadid2020,Cueto2009},  etc. 

Nonlinearity and discretization of high-order spatial operators (fourth-order) are the main issues encountered when numerically solving these equations. The nonlinearity is often induced by non-constant mobility and nonlinear free energy functions. As a result, several different procedures have been employed over the years to deal with higher-order equations, for instance, mixed formulations, where the governing equations are split into a coupled system of lower-order differential equations \cite{BoffiBrezziFortin2013,Elliott1989}.

There have been many studies on numerical methods for solving parabolic equations involving the bi-harmonic operator, including finite difference \cite{Bertozzi1,SHIN2011,Bertozzi3}, finite volume \cite{DARGAVILLE2015,Grun2000}, Fourier spectral \cite{HE2009,YE2005} and finite element methods \cite{Barrett1,Brenner2012,Dziuk2008,JOKISAARI2017}. Most studies using finite element method focus on mixed formulation to approximate the CH equations \cite{Barrett3,Barrett2,DU2011,Goudenege2012,KAESSMAIR2016,KASTNER2016,ZHANG2013}, see also the references therein. Second-order time-accurate methods for the CH equations have mostly been studied for variants of the Crank-Nicolson scheme with various treatment of the nonlinear energy term \cite{Tierra2013,Qiang1991,Vignal2017,Wu2014}. Convex-concave splitting of the energy term is often employed, but it leads to nonlinear schemes except for some particular forms of the energy function. Addition of stabilization terms, use of Hermite quadrature formula and/or explicit treatment or extrapolation of the energy function, although leading to a linear scheme, is very limited in term of stability constraint \cite{Tierra2013,Shen2010}. Recently, Yang et al.\ \cite{Yang2019} proposed a family of second-order unconditionally energy-stable schemes for CH equations by reformulating the equations into an equivalent system employing a scalar auxiliary variable, giving rise  to four equations to solve at each time step.

In the current study, we focus on the numerical approximation of the CH equations. We propose a novel linear semi-implicit time stepping scheme to approximate the  problem. The proposed scheme maintains second-order accuracy in time by using second-order backward difference formula (BDF2) for the time derivative, extrapolation formula for the non-constant mobility function and a special technique using semi-implicit Taylor approximation for the nonlinear energy function. The use of Taylor expansion on a term in solving partial differential equations was introduced by Rosenbrock \cite{Rosenbrock1963}. A scheme employing the Crank-Nicolson time integration and Taylor series expansions of the nonlinear terms were proposed by Vignal et al.\ \cite{Vignal2017} for general classes of phase-field models with polynomial potentials.  Their resulting system is nonlinear and therefore, requires an iterative method for the linearization. Our scheme applies to general classes of phase-field models with any form of free energy function of class $\mathcal{C}^2$ and does not require an iterative method. A linear system is solved at each time step. We perform numerical convergence analysis in space using both a manufactured and analytical solutions to demonstrate optimal convergence rates of the $L^2$- and $H^1$-error norms using standard Lagrange elements. Better, super-convergence of the $L^2$-error norm is observed with quadratic Lagrange elements. Numerical results from selected test cases also show that the proposed scheme is efficient in terms of computational cost and  satisfies the free energy dissipation and the mass conservation properties.

The outline of the paper is as follows.  In section~{\ref{SectModelEquation}}, the model equation is introduced and its main mathematical properties are briefly discussed. In section~{\ref{SectNumMeth}}, the fully discrete finite element approximation scheme is presented. In Section~{\ref{SectNumResults}}, the spatial convergence rates, accuracy and robustness of the proposed scheme are numerically investigated by means of several numerical test cases. Second-order time-accuracy is verified with both constant and variable mobility function, and the discrete energy dissipation and the discrete mass conservation properties of the scheme are illustrated. Section~{\ref{SectConclusion}} provides a summary and some concluding remarks.

\section{Model equation}
\label{SectModelEquation}
Let $\Omega\subset\mathbb{R}^d$, where $d=2$ or $3$, be an open bounded set with sufficiently smooth boundary $\partial \Omega$, and $T>0$ a fixed positive time. We focus on a class of fourth-order parabolic PDEs of the form
\begin{equation}
\partial_tu=-\gamma\nabla\cdot\big(f(u)\nabla \Delta u\big) + \nabla\cdot\big(f(u)\nabla \varphi'(u)\big) \quad \text{in } \Omega\times(0,T),
\label{ThinFilmEquation}
\end{equation}
with unknown $u:\Omega\times[0,T]\rightarrow\mathbb{R}$, the given functions $f:\mathbb{R}\rightarrow\mathbb{R}_{\geqslant 0}$, $\varphi:\mathbb{R}\rightarrow\mathbb{R}$ are  smooth, and $\gamma\in\mathbb{R}_{\geqslant 0}$ is a constant. Depending on the applications, $u$ can describe the height of a liquid phase spreading on a solid surface \cite{greenspan1978,Hocking1981}, the concentration, volume fraction or density of a phase in a mixture \cite{CAHN1962,CahnHilliard1958}, etc.

By introducing an auxiliary unknown  $w$, the fourth-order equation \eqref{ThinFilmEquation} can be written in the following system of second-order equations
\begin{equation}
\begin{aligned}
\partial_tu&=\nabla\cdot\big(f(u)\nabla w\big) \quad \text{in } \Omega\times(0,T),\\
w&=-\gamma\Delta u +\varphi'(u) \quad \text{in } \Omega\times(0,T).
\end{aligned}
\label{ThinFilmSyst}
\end{equation}
System \eqref{ThinFilmSyst} is to be supplemented by the initial condition
\begin{equation}
u(\bm{x},0)=u_0(\bm{x}) \quad \text{in } \Omega,
\label{ThinFilmSystInitialCondition}
\end{equation}
and the no-flux boundary conditions
\begin{equation}
\nabla u\cdot\bm{n}=\nabla w\cdot\bm{n}=0 \quad \text{on } \partial\Omega\times(0,T),\\
\label{HomoNeumannBoundCond}
\end{equation}
where $\bm{n}$ is the exterior unit normal vector to $\partial\Omega$.

There are two important properties for the  equations \eqref{ThinFilmSyst} with the natural boundary conditions \eqref{HomoNeumannBoundCond}. The first one is the conservation of mass. In fact,
\begin{equation}
\dfrac{d}{dt}\int_{\Omega}u\,d\bm{x}=\int_{\Omega}\partial_tu\,d\bm{x}=\int_{\Omega}\nabla\cdot\big(f(u)\nabla w\big)\,d\bm{x}=\int_{\partial\Omega}f(u)\nabla w\cdot\bm{n}\,d\bm{x}=0,
\end{equation}
that is the total amount of the phase in the domain must always be equal to the given initial amount of this phase. The equation is associated with the total free energy functional
\begin{equation}
J(u)=\int_{\Omega}\left(\frac{\gamma}{2}\vert\nabla u\vert^2+\varphi(u)\right)\,d\bm{x},
\label{EnergyFunctional}
\end{equation}
which satisfies the total free energy dissipation, that is
\begin{equation}
\dfrac{d J}{dt}=\int_{\Omega}\big(\gamma\nabla u \nabla \partial_tu+\varphi'(u)\partial_tu\big)\,d\bm{x}=\int_{\Omega}w \partial_tu\,d\bm{x}=-\int_{\Omega}f(u)\vert\nabla w\vert^2\,d\bm{x} \leqslant 0.
\label{EnergyProperty}
\end{equation}

A weak formulation of problem \eqref{ThinFilmEquation} can be obtained by multiplying \eqref{ThinFilmEquation} by test function $v$. Integrating by parts and using the boundary conditions \eqref{HomoNeumannBoundCond} eventually lead to: Find $u \in H^2(\Omega)$ such that
\begin{equation}
\int_{\Omega}\partial_tu v\,d\bm{x}+\int_{\Omega}\gamma f(u)\Delta u\Delta v\,d\bm{x}+\int_{\Omega}\gamma \Delta u\nabla f(u)\cdot\nabla v\,d\bm{x}+\int_{\Omega}f(u)\nabla \varphi'(u)\cdot\nabla v\,d\bm{x}=0,
\label{ThinEquationWeakForm}
\end{equation}
for all $v \in H^2(\Omega)$. A weak formulation for problem \eqref{ThinFilmSyst} can be obtained by multiplying the first and second equation in \eqref{ThinFilmSyst} by test functions $v$ and $q$, respectively: Find $\big(u ,w\big) \in H^1(\Omega)\times H^1(\Omega)$ such that
\begin{equation}
\begin{aligned}
&\int_{\Omega}u_t v \,d\bm{x}+\int_{\Omega} f(u)\nabla w\cdot\nabla v\, d\bm{x}=0,\quad\forall v \in H^1(\Omega),\\
&\int_{\Omega} w q\,d\bm{x}-\int_{\Omega}\gamma\nabla u\cdot\nabla q\,d\bm{x}-\int_{\Omega}\varphi'(u) q\,d\bm{x}=0,\quad\forall q \in  H^1(\Omega).
\end{aligned}
\label{lubweakForm}
\end{equation}
Theoretical results on the existence of solutions to \eqref{ThinFilmEquation} and \eqref{ThinFilmSyst} in the sense \eqref{ThinEquationWeakForm} and \eqref{lubweakForm}, respectively, can be found in \cite{Dai2016,Elliott89,Elliott96}.

\section{Numerical methods}
\label{SectNumMeth}
\subsection{Spatial discretization}
K\"{a}stner et al.\ \cite{KASTNER2016} showed that the convergence rates for approximations that are less than second-order in space are suboptimal when using the direct weak formulation \eqref{ThinEquationWeakForm}. This is due to the discretization of higher-order differential operators.  The weak formulation \eqref{lubweakForm} naturally leads to mixed finite element methods which provide robust numerical approaches to compute solutions of high-order PDEs by avoiding the discretization of high-order derivatives \cite{BoffiBrezziFortin2013,Elliott1989}.  In the following, we will use the mixed weak formulation \eqref{lubweakForm}.

We discretize in space by continuous piecewise Lagrange finite elements. Given a polygonal domain $\Omega$, we consider  $\mathcal{T}_h$ being a uniform partitioning of the domain $\Omega$ into disjoint triangles $\kappa$, so that 
\begin{equation}
\overline{\Omega}=\bigcup_{\kappa\in\mathcal{T}_h}\kappa.
\end{equation}
The  element mesh size $h$ is defined as 
\begin{equation}
h:=\max_{\kappa\in \mathcal{T}_h}h_{\kappa},
\end{equation}
where $h_{\kappa}:=diam(\kappa)$. Associated with the partitioning $\mathcal{T}_h$  is the finite element space
\begin{align}
&{\mathcal{V}_h}=\Big\lbrace v_h\in \mathcal{C}^0(\overline{\Omega},\mathbb{R}): {v_h}_{|_{\kappa}}\in \mathbb{P}_k,\quad \forall \kappa\in \mathcal{T}_h \Big\rbrace\subset H^1(\Omega),
\end{align}
where $\mathbb{P}_k$ denotes the space of polynomials of degree less than or equal to $k$ on any element $\kappa$. 

The resulting semi-discrete (in space) weak formulation of problem \eqref{ThinFilmSyst}-\eqref{HomoNeumannBoundCond} is given as follows: Find $\big( u_h,w_h\big)\in \mathcal{C}^1\big(0,T;{\mathcal{V}_h} \times {\mathcal{V}_h}\big)$ such that, for all $t\in[0,T]$,
\begin{equation}
\begin{aligned}
&\int_{\Omega}\partial_tu_h(t)v_h \,d\bm{x}+\int_{\Omega} f(u_h(t))\nabla w_h(t)\cdot\nabla v_h\,d\bm{x}=0,\quad\forall v_h  \in {\mathcal{V}_h} ,\\
&\int_{\Omega}w_h(t) q_h\,d\bm{x}-\int_{\Omega}\gamma\nabla u_h(t)\cdot\nabla q_h\,d\bm{x}-\int_{\Omega}\varphi'(u_h(t)) q_h\,d\bm{x}=0,\quad\forall q_h \in {\mathcal{V}_h},\\
&u_h(0)=\Pi_hu_0,
\end{aligned}
\label{ThinSystSemidiscreteWeakForm}
\end{equation}
where $\Pi_h$ is an interpolation or projection operator on ${\mathcal{V}_h}$.

\subsection{Temporal discretization}
The time interval $[0,T]$ is discretized as
\begin{equation}
t_n=n \Delta t,\quad  n=0,1,...,N,\qquad  \Delta t=\dfrac{T}{N}, 
\end{equation} 
where $\Delta t$ is the time step used. We consider a uniform discretization in time by a semi-implicit stepping method, so that $u_h^n\simeq u_h(t_n)$ for $n=0,1,..., N$. For the temporal discretization of the equations \eqref{ThinFilmSyst}, all linear terms  are approximated implicitly. The time derivative $\partial_tu_h$ is approximated using a second-order backward differentiation formula (BDF2). To allow for a linear scheme, a second-order extrapolation formula is employed for the mobility function $f(u_h)$ while the energy term $\varphi'(u_h)$ is evaluated using the Taylor approximation of second-order truncation error. We have noticed from numerical experiments that the approximation of the nonlinear energy term using the proposed technique is better than by an extrapolation which leads to severe stability constraint. Our technique introduces a linear implicit term which improves the stability constraint. We emphasize two main benefits resulting from the proposed scheme. First of all, it is easy to implement. No fixed point iteration for solving the system is required. Only one linear system needs to be solved at every time step. The proposed scheme is written as follows: given a suitable approximation of the initial solution $u_h^{-1}=\Pi_hu_0\in {\mathcal{V}}_h$ and a proper initialization for $u_h^0\in {\mathcal{V}}_h$, find $\big(u_h^n,w_h^n\big)\in {\mathcal{V}}_h \times {\mathcal{V}}_h$ such that
\begin{eqnarray}
\begin{aligned}
&\int_{\Omega}\left(\dfrac{3u_h^n-4u_h^{n-1}+u_h^{n-2}}{2\Delta t}\right)v_h\,d\bm{x}+\int_{\Omega}\big(2f(u_h^{n-1})-f(u_h^{n-2})\big)\nabla w_h^n\cdot\nabla v_h\,d\bm{x}=0,\quad\forall v_h \in {\mathcal{V}_h},\\
&\int_{\Omega} w_h^nq_h \,d\bm{x}-\int_{\Omega}\gamma\nabla u_h^n\cdot\nabla q_h \,d\bm{x}-\int_{\Omega}\Big(\varphi'(u_h^{n-1})+\varphi''(u_h^{n-1})\big(u_h^n-u_h^{n-1}\big)\Big)q_h\,d\bm{x}=0, \quad\forall q_h \in {\mathcal{V}_h},
\end{aligned}
\label{SBDF2ThinWeakForm}
\end{eqnarray}
for all $1\leqslant n\leqslant N$. \\
The numerical algorithm \eqref{SBDF2ThinWeakForm} is a two steps scheme and therefore, it requires the use of a starting procedure to obtain $u_h^0$. Here, we use a semi-implicit backward Euler method 
\begin{eqnarray}
\begin{aligned}
&\int_{\Omega}\left(\dfrac{u_h^0-u_h^{-1}}{\Delta t}\right)v_h\,d\bm{x}+\int_{\Omega}f(u_h^{-1})\nabla w_h^0\cdot\nabla v_h\,d\bm{x}=0,\quad\forall v_h \in {\mathcal{V}_h},\\
&\int_{\Omega} w_h^0q_h \,d\bm{x}-\int_{\Omega}\gamma\nabla u_h^0\cdot\nabla q_h \,d\bm{x}-\int_{\Omega}\varphi'(u_h^{-1})q_h\,d\bm{x}=0, \quad\forall q_h \in {\mathcal{V}_h},
\end{aligned}
\label{SBDF1ThinWeakForm}
\end{eqnarray}
to compute $u_h^0$ starting from  $u_h^{-1}$.

The numerical treatment of the energy term $\varphi'(u_h)$ uses the following semi-implicit Taylor expansion of $\varphi'(u_h(t_n))$ about $u_h=u_h(t_{n-1})$: 
\begin{equation}
\varphi'(u_h(t_{n-1}))+\varphi''(u_h(t_{n-1}))\big(u_h(t_n)-u_h(t_{n-1})\big)=\varphi'(u_h(t_n))-\frac{1}{2}\big(u(t_n)-u(t_{n-1})\big)^2\varphi^{(3)}(\overline{u}^n),
\label{DirectExpansion}
\end{equation}
where $\overline{u}^n$ is between $u_h(t_{n-1})$ and $u_h(t_n)$. To evaluate the local truncation error
\begin{equation}
\tau^n(q_h)=\int_{\Omega}-\frac{1}{2}\big(u(t_n)-u(t_{n-1})\big)^2\varphi^{(3)}(\overline{u}^n)q_h\,d\bm{x}
\label{TruncErrorEnergy}
\end{equation}
of the energy term, one expands again $u_h(t_{n-1})$ in \eqref{TruncErrorEnergy} around the point $t=t_n$ and gets
\begin{equation}
\tau^n(q_h)=-\frac{\Delta t^2}{2}\int_{\Omega}\big(\partial_tu_h(\xi^n)\big)^2\varphi^{(3)}(\overline{u}^n)q_h\,d\bm{x},
\label{tau3final}
\end{equation}
where $\xi^n=\xi^n(\bm{x})\in[t_{n-1},t_n]$ for all $\bm{x}\in\Omega$. The backward differentiation formula used for the numerical treatment of $\partial_t u_h$ is known to have error that is second-order in time and one can easily show that the extrapolation technique used for the approximation of the mobility term $f(u_h)$ has an order of accuracy of $\mathcal{O}(\Delta t^2)$. Hence, the scheme \eqref{SBDF2ThinWeakForm} is second-order accurate in time.

\subsection{Solution of the discrete system}
Equations in \eqref{SBDF2ThinWeakForm} or \eqref{SBDF1ThinWeakForm} consist in a coupled linear system in $u_h^n$ and $w_h^n$ of the following form
\begin{equation}
\begin{bmatrix}
M_1 & A_1 \\
-A_2 & M_2
\end{bmatrix}
\begin{bmatrix}
U \\
W
\end{bmatrix}
=
\begin{bmatrix}
F \\
G
\end{bmatrix},
\end{equation}
where $M_i$ and $A_i$, $i=1,2$, are mass and stiffness matrices, respectively; $U$ and $W$ correspond to the finite element degrees of freedom for $u_h^n$ and $w_h^n$, respectively; $F$ and $G$ are properly chosen constant vectors.

This linear system shares characteristics with those for saddle point problems, such as the block-structured sparsity pattern, the lack of positive definiteness, etc. These linear systems are difficult to solve. The use of a direct solver may lead to a lot of fill-in, while iterative solvers may be impacted by the indefiniteness of the matrix. For our 2D test cases, we used the sparse linear solver UMFPACK \cite{Davis2004}, which mitigates the fill-in and provides a relatively efficient solution. However, for 3D computations, one would need to use iterative methods. One may attempt to build a dedicated iterative method accounting for the fact that all the blocks $M_i$ and $A_i$ are at least semi-positive definite, a property usually not shared by general saddle point problems.

\subsection{Mass conservation}
As usual in finite element methods with no-flux boundary conditions, the scheme \eqref{SBDF2ThinWeakForm} is conservative. In fact, taking $v_h=1$ in the first equation of \eqref{SBDF2ThinWeakForm} and \eqref{SBDF1ThinWeakForm}, respectively, leads to
\begin{equation}
\int_{\Omega}\dfrac{3u_h^n-4u_h^{n-1}+u_h^{n-2}}{2\Delta t}\,d\bm{x}=0 \qquad\text{and}\qquad \int_{\Omega}\dfrac{u_h^0-u_h^{-1}}{\Delta t}\,d\bm{x}=0.
\label{SchemeConservation}
\end{equation}
It follows by induction from \eqref{SchemeConservation} that 
\begin{equation}
\int_{\Omega}u_h^n\,d\bm{x}=\int_{\Omega}u_h^0\,d\bm{x},\quad \forall n=1,...,N,
\label{SchemeConservation5}
\end{equation}
which means that the discrete total mass is conserved over time.

\section{Numerical results}
\label{SectNumResults}
A series of numerical problems in two dimensions is given in this section to test the performance of the numerical method developed in section \ref{SectNumMeth}. All the algorithms for the simulations are implemented using FreeFem++ software \cite{FreeFem}.

\subsection{Numerical results --- Convergence}
In this subsection we numerically show the rates of convergence in space and time of the algorithm \eqref{SBDF2ThinWeakForm}-\eqref{SBDF1ThinWeakForm}. We first conduct a detailed investigation using the method of manufactured solution for a constant mobility function. We then use an analytic solution for a non-constant mobility function to perform a similar error analysis.

\subsubsection{Manufactured solution}
\label{SubsectionManufSol}
Following the framework presented in \cite{KAESSMAIR2016,KASTNER2016,Yang2019,ZHANG2013}, we employ a manufactured solution to test our method. The idea is to use an arbitrary regular function as analytic solution for the computation of numerical errors.

Substituting in \eqref{ThinFilmSyst} the  function
\begin{equation}
\widehat{u}=(t+1)\sin(\alpha\pi x),
\label{ManufSolu}
\end{equation}
where $\alpha\in\mathbb{R}$, by assuming a free energy function
\begin{equation}
\varphi(u)=\frac{u^4}{4}-\frac{u^2}{2},
\label{FreeEnergyFunction}
\end{equation}
and a constant mobility function $f(u)=1$, yields 
\begin{equation}
\left\lbrace
\begin{aligned}
&\widehat{u}_t=\nabla\cdot(\nabla\widehat{w})+S,\\
&\widehat{w}=-\gamma\Delta\widehat{u}+\varphi'(\widehat{u}),
\end{aligned}
\right.
\label{ManufacSollubSyst}
\end{equation}
where the residual $S$ is given by
\begin{equation}
\begin{aligned}
S=&\sin(\alpha\pi x) + 3(t+1)^3\alpha^2\pi^2\sin^3(\alpha\pi x) + \gamma(t+1)\alpha^4\pi^4\sin(\alpha\pi x) - (t+1)\alpha^2\pi^2\sin(\alpha\pi x)\\
&-6(t+1)^3\alpha^2\pi^2\sin(\alpha\pi x)\cos^2(\alpha\pi x)
\end{aligned}
\label{Residual}
\end{equation}
and the auxiliary variable $\widehat{w}$ is given by
\begin{equation}
\widehat{w}=\gamma(t+1)\alpha^2\pi^2\sin(\alpha\pi x) + (t+1)^3\sin^3(\alpha\pi x) - (t+1)\sin(\alpha\pi x).
\label{ManufSolw}
\end{equation}
We assumed here a constant mobility function $f(u)\equiv 1$ to simplify the calculation of the residual function $S$, but the energy term \eqref{FreeEnergyFunction} is a nonlinear non-zero function. This will allow to test the impact of the new proposed technique for the discretization of the nonlinear energy term, based on a semi-implicit Taylor expansion of the latter.

For the convergence analysis in space, we compute the $L^2$- and $H^1$-error norms  on the discrete solution $\widehat{u}_h$ of \eqref{ManufacSollubSyst} at the final time, which are defined as 
\begin{equation}
\begin{aligned}
&\Vert\widehat{u}-\widehat{u}_h\Vert_{L^2(\Omega)}=\left(\int_{\Omega}\vert\widehat{u}-\widehat{u}_h\vert^2\,d\bm{x}\right)^{1/2},\\
&\Vert\widehat{u}-\widehat{u}_h\Vert_{H^1(\Omega)}=\left(\int_{\Omega}\Big(\vert\widehat{u}-\widehat{u}_h\vert^2+\vert\nabla(\widehat{u}-\widehat{u}_h)\vert^2\Big)\,d\bm{x}\right)^{1/2}.
\end{aligned}
\end{equation}
A set of boundary conditions for the system has to be imposed. We will investigate both Dirichlet and non-homogeneous Neumann boundary conditions using the exact solution \eqref{ManufSolu} and \eqref{ManufSolw} at boundaries. The linear ($\mathbb{P}_1$) and quadratic ($\mathbb{P}_2$) Lagrange finite element discretizations and four uniform meshes ($25\times 25$, $50\times 50$, $100\times 100$ and $200\times 200$ elements) of the unit square (centered at the origin) domain  have been considered. We use a fixed time step $\Delta t= 10^{-5}$ and consider a small value for the parameter $\gamma=10^{-4}$. The final computational time is $T=10^{-3}$.

\begin{figure}[!h]
\centering
\begin{tabular}{cc}
\includegraphics[scale=0.57]{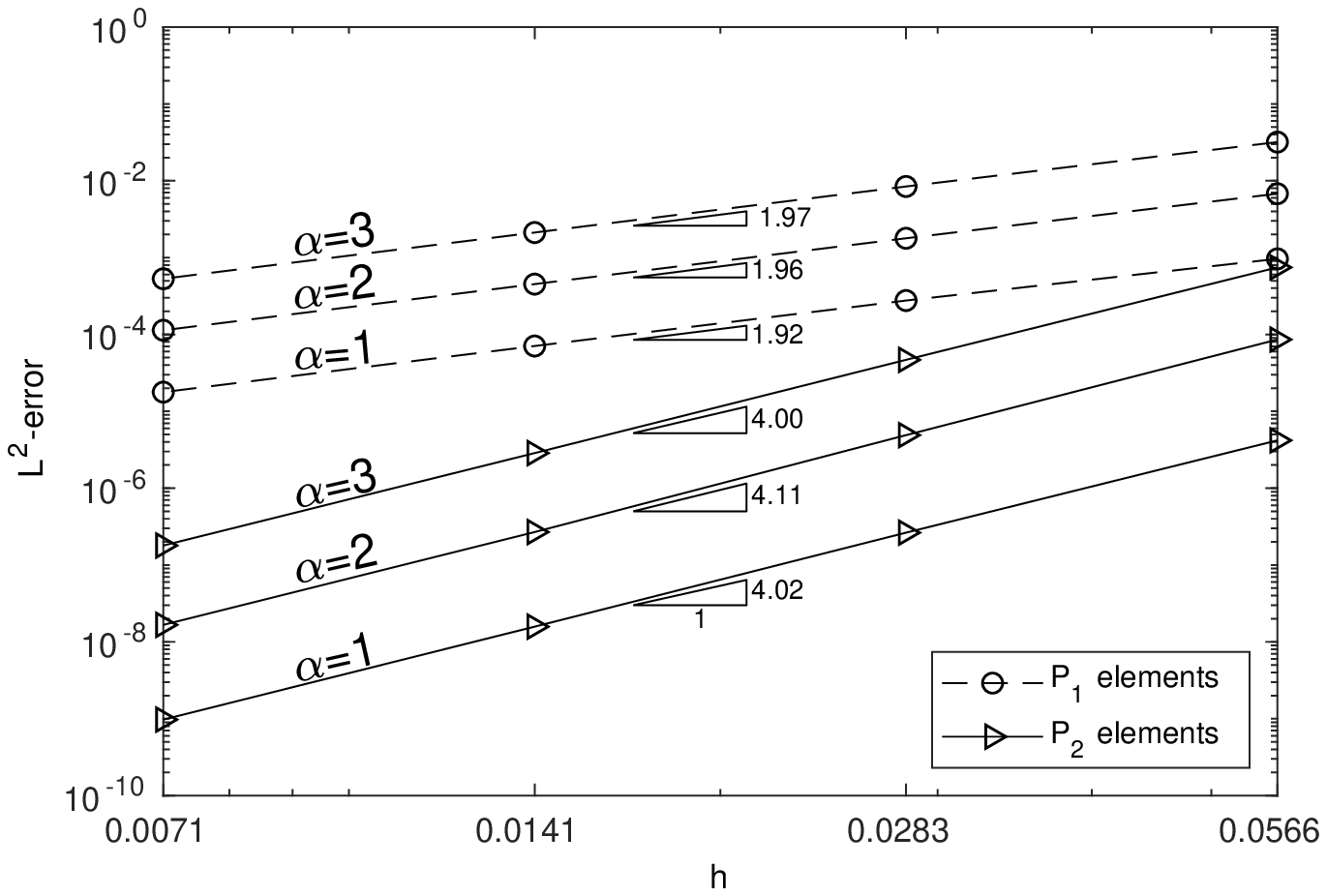}&\includegraphics[scale=0.57]{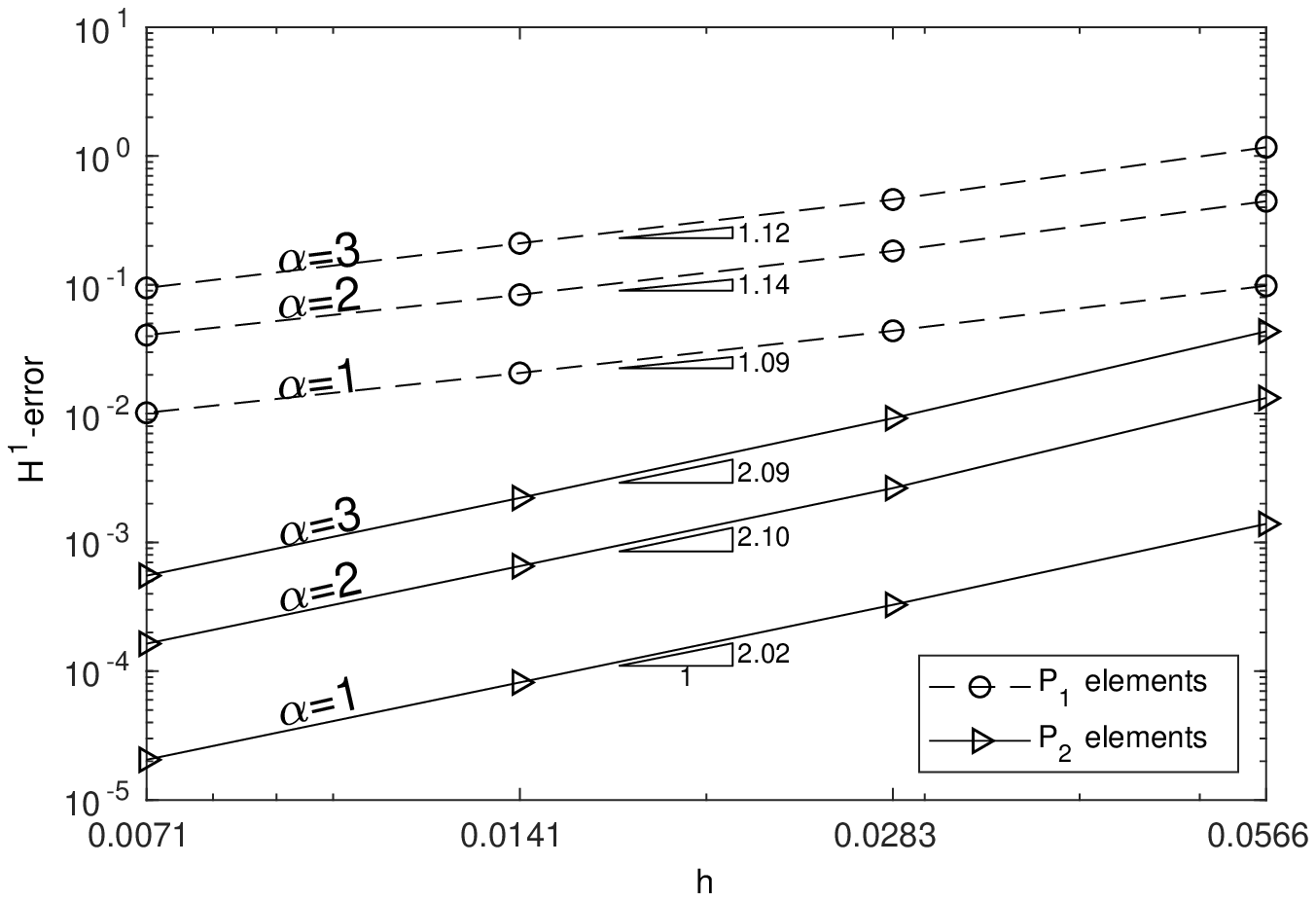}\\
$(a)$ & $(b)$
\end{tabular}
\caption{Error and convergence rates of linear and quadratic Lagrange elements: $(a)$ $L^2$-error, and $(b)$ $H^1$-error on the discrete solution $\widehat{u}_h$ as function of the element size $h$. Results are obtained with Dirichlet boundary conditions.}
\label{Fig:L2errorH1errorP1P2elementsDirichlet}
\end{figure}

We first impose Dirichlet conditions $\widehat{u}_h=\Pi_h\widehat{u}$ and $\widehat{w}_h=\Pi_h\widehat{w}$ on boundaries.  The numerical errors on $\widehat{u}_h$ are computed for different values of the parameter $\alpha$. Figure~\ref{Fig:L2errorH1errorP1P2elementsDirichlet} shows theses errors as a function of the element size $h$. Optimal convergence rates, i.e.\ $\mathcal{O}(h^{k+1})$ and $\mathcal{O}(h^k)$,  are obtained for the $L^2$- and $H^1$-error, respectively. Better, we get super convergence for the $L^2$-norm error with the quadratic Lagrange elements. The error increases with $\alpha$, which controls the gradient of the solution. However, the rates of convergence remain the same, and therefore, are independent from the choice of $\alpha$. 

\begin{figure}[!h]
\centering
\begin{tabular}{cc}
\includegraphics[scale=0.57]{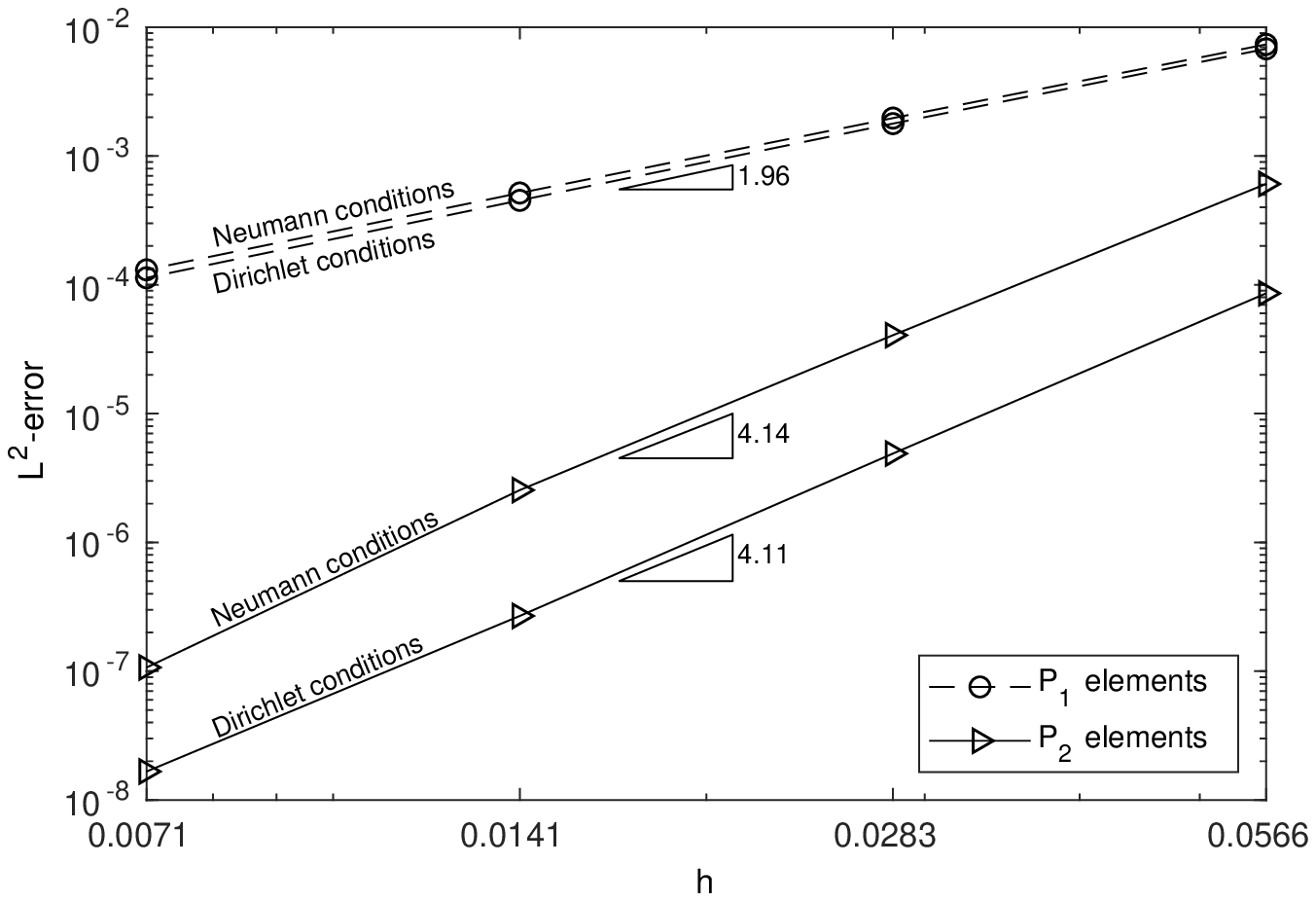}&\includegraphics[scale=0.57]{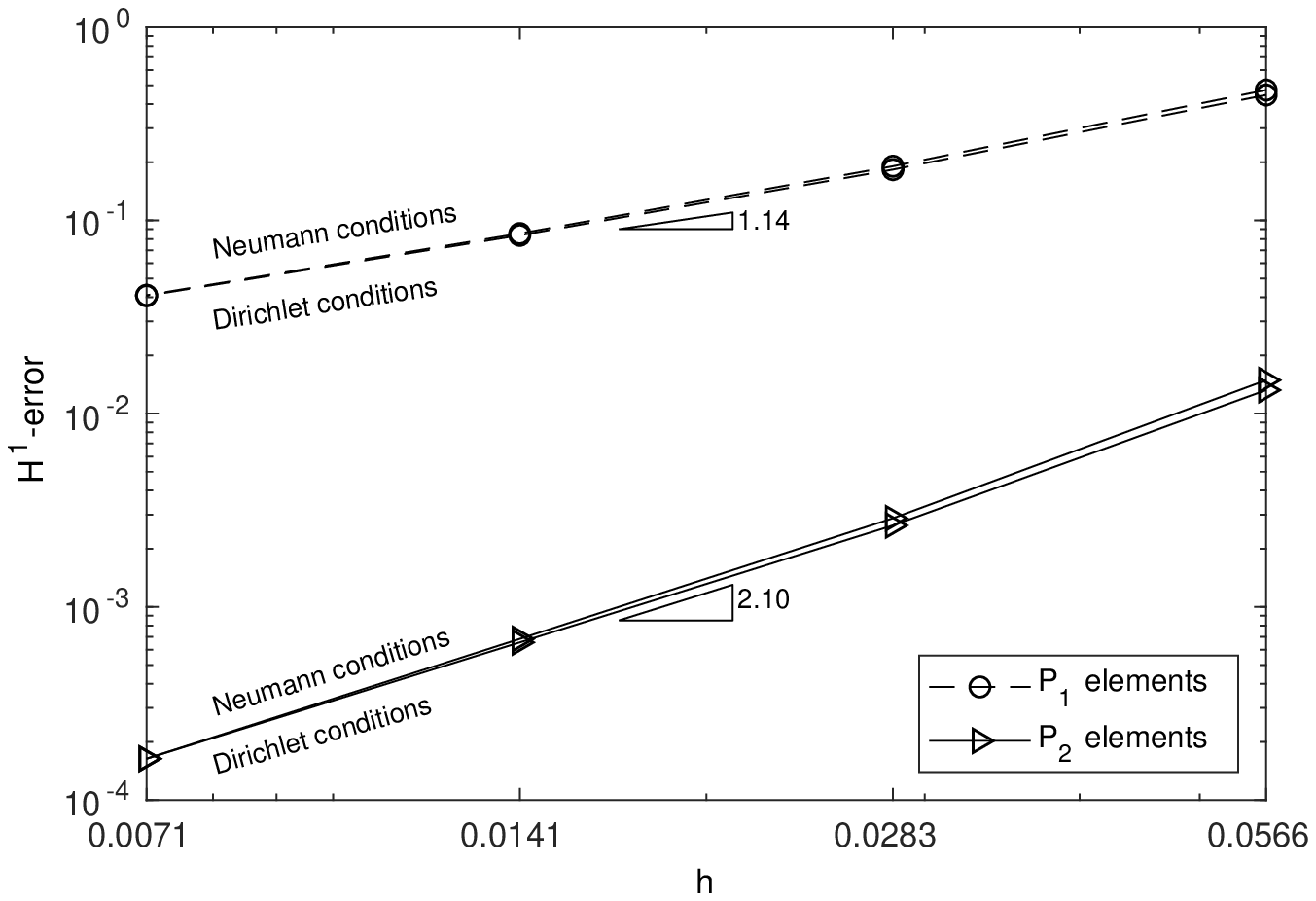}\\
$(a)$&$(b)$
\end{tabular}
\caption{Quantitative comparison of the error computed with Dirichlet and Neumann boundary conditions, respectively: $(a)$ $L^2$-error, and $(b)$ $H^1$-error on the discrete solution $\widehat{u}_h$ as function of the element size $h$.}
\label{Fig:L2errorH1errorP1P2elementsDirichletVsNeumann}
\end{figure}

We next impose non-homogeneous Neumann conditions $\frac{\partial \widehat{u}_h}{\partial\mathbf{n}}=\Pi_h\frac{\partial \widehat{u}}{\partial\mathbf{n}}$ and $\frac{\partial \widehat{w}_h}{\partial\mathbf{n}}=\Pi_h\frac{\partial \widehat{w}}{\partial\mathbf{n}}$. Here, also optimal convergence rates in the $L^2$- and $H^1$-error are obtained. Again, super-convergence in the $L^2$-norm is observed with quadratic Lagrange elements. To facilitate comparison between the results obtained using Dirichlet and Neumann conditions, the errors computed for $\alpha=2$ are shown in Figure~\ref{Fig:L2errorH1errorP1P2elementsDirichletVsNeumann}.  We can see that the $L^2$-error with Neumann conditions are larger than with Dirichlet conditions, but the $H^1$-error are almost identical.

We have performed a similar error analysis with the $L^2$- and $H^1$-error on the auxiliary variable $\widehat{w}_h$ and found the same orders of convergence as for $\widehat{u}_h$. However, the $H^1$-error on $w_h$ is larger and the order of magnitude of the difference increases with $\alpha$. This is natural since $\widehat{w}_h=-\gamma\Delta \widehat{u}_h+\varphi'(\widehat{u}_h)$ is more difficult to approximate and the gradient of the solution increases with $\alpha$ more rapidly than for $\widehat{u}_h$.

\begin{figure}[!h]
\centering
\begin{tabular}{cc}
\includegraphics[scale=0.57]{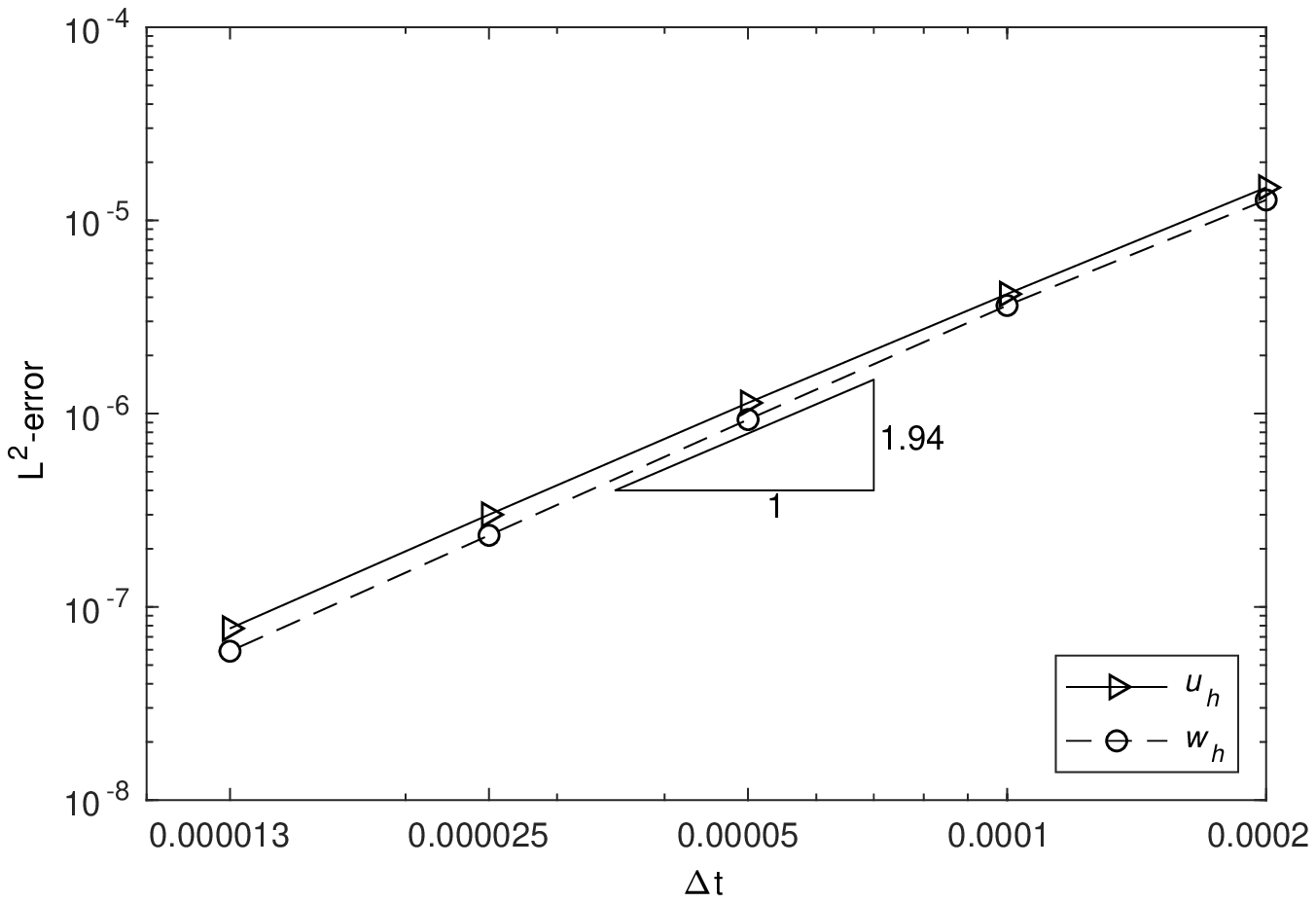}&\includegraphics[scale=0.57]{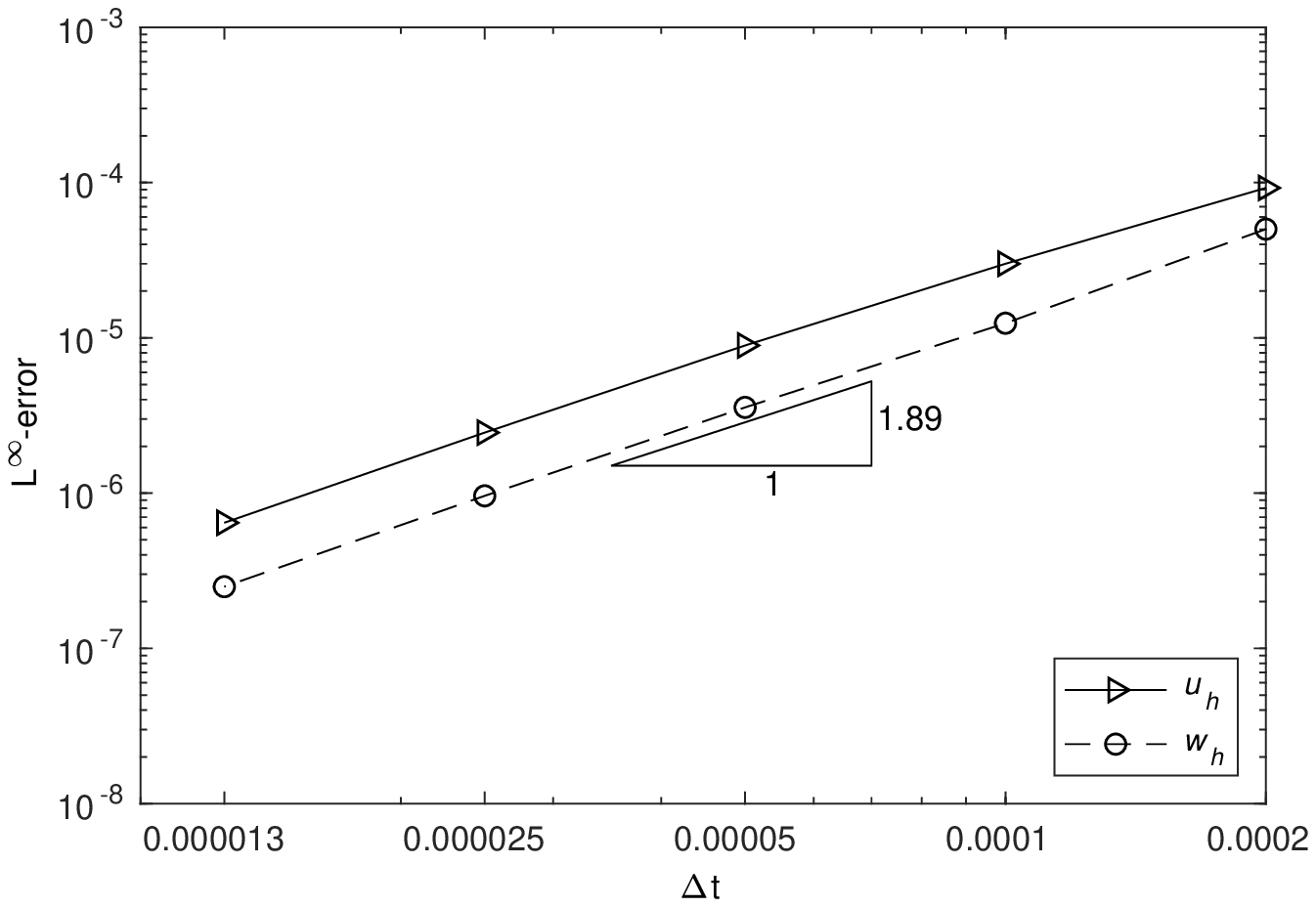}\\
$(a)$& $(b)$
\end{tabular}
\caption{Temporal convergence rate: $(a)$ $L^2$-error, and $(b)$ $L^{\infty}$-error on $\widehat{u}_h$ and $\widehat{w}_h$, respectively, as function of time step size $\Delta t$. Results are obtained with $\mathbb{P}_2$ elements and Dirichlet conditions.}
\label{Fig:L2errorLinftyerrorManufSolTimeOrderSBDF2}
\end{figure}

To test the temporal convergence rate, we compute a reference solution ($u_h^{\text{ref}}$,$w_h^{\text{ref}}$) on a grid of $50\times 50$ elements using a very small time step $\Delta t= 10^{-8}$ at the final time $T=10^{-3}$. We choose the quadratic Lagrange element and impose Dirichlet boundary conditions to minimize the spatial error. Moreover, with quadratic elements, fixed computational time and time step give smaller error than with $\mathbb{P}_1$ elements (see Table~\ref{tab:ExactSolL2H1errorP1P2Dirichlet} in section~\ref{SubsectionAnalyticSol}).  We analyze the numerical convergence in time by varying the time step $\Delta t$ systematically between $0.125\times 10^{-4}$ and $2\times 10^{-4}$. For each time step, the numerical errors in the $L^2$- and $L^{\infty}$-norms  are computed. Figure~\ref{Fig:L2errorLinftyerrorManufSolTimeOrderSBDF2} shows the $L^2$- and $L^{\infty}$-error on $\widehat{u}_h$ and $\widehat{w}_h$, respectively, and the order of convergence of the numerical solution with respect to the reference solution.  We obtain almost a second-order convergence rate in time.

\subsubsection{Analytic solution}
\label{SubsectionAnalyticSol}
In the previous subsection, optimal convergence rates have been obtained for constant mobility function. It is the main purpose of this subsection to numerically study the convergence rates for non-constant mobility function by using a particular analytic solution.

Setting $f(u)=u$, $\varphi=0$ and $\gamma=1$ in \eqref{ThinFilmSyst}, we get the nonlinear fourth-order equation
\begin{equation}
u_t+\nabla\cdot(u\nabla\Delta u)=0 \quad \text{in } \Omega\times(0,T).\\
\label{lubSyst}
\end{equation} 
Equation \eqref{lubSyst} has a nonnegative closed-form compactly-supported, $2$-dimensional self-similar solution \cite{ferreira_bernis_1997} given by
\begin{equation}
u(x,y,t)=
\left\lbrace
\begin{aligned}
&\dfrac{t^{-1/3}}{192}(L^2-r^2)^2,\quad r<L,\\
&0,\quad r\geqslant L,
\end{aligned}
\right.
\qquad\quad r=\dfrac{\sqrt{x^2+y^2}}{t^{1/6}}.
\label{ExactSol1}
\end{equation}
In the framework of the mixed formulation, we need the auxiliary variable $w=-\Delta u$, which is calculated in a straightforward way from \eqref{ExactSol1} to get
\begin{equation}
w(x,y,t)=
\left\lbrace
\begin{aligned}
&\dfrac{t^{-2/3}}{24}(L^2-2r^2),\quad r<L,\\
&0,\quad r\geqslant L.
\end{aligned}
\right.
\label{ExactSol2}
\end{equation}
The auxiliary variable \eqref{ExactSol2} is not continuous at point where $r=L$. Moreover, the source-type solution \eqref{ExactSol1} tends to a Dirac mass as $t$ goes to zero. Therefore, we have to choose $L$ large enough and  a starting time $t_0>0$ to ensure the regularity of the solution \eqref{ExactSol1}-\eqref{ExactSol2} for the mixed formulation. Here, there is no energy term, but the mobility function depends on the solution. This will allow to test the impact of the numerical treatment for the  mobility function, based on extrapolation. 

To test the spatial convergence rates, we compute the $L^2$- and $H^1$-error norms between the exact solution \eqref{ExactSol1}--\eqref{ExactSol2} and the numerical solution computed with the proposed method \eqref{SBDF2ThinWeakForm}-\eqref{SBDF1ThinWeakForm}, where we set $\varphi\equiv0$. We use the same computational domain, meshes, finite element approximations as in section \ref{SubsectionManufSol}. Since, we focus on the spatial convergence rates, we use a small time step $\Delta t=10^{-8}$ to minimize the temporal error from the extrapolation of the non constant mobility. We choose $L=3$ and a starting time $t_0=0.001$. Tables~\ref{tab:ExactSolL2H1errorP1P2Dirichlet} and \ref{tab:ExactSolL2H1errorP1P2Neumann} show  the  computational time, the $L^2$- and $H^1$-error as function of the element $h$, computed with Dirichlet and Neumann boundary conditions, respectively. The formula to obtain the order of convergence in space, $p$, reads:
\begin{equation}
p=\ln\left(\dfrac{\Vert u-u_h\Vert}{\Vert u-u_{h/\tau}\Vert}\right) \biggm/ \ln \tau,
\label{ConvergenceRateFormula}
\end{equation}
where $\tau>0$ is the refinement factor between consecutive element sizes. The value $\tau=2$ was chosen in our numerical experiments. By using the formula \eqref{ConvergenceRateFormula}, we  get from Table \ref{tab:ExactSolL2H1errorP1P2Dirichlet} that the convergence rates in the $H^1$-norm are about $1$ and $2$ for the linear and quadratic Lagrange elements, respectively, while the convergence rates in the $L^2$-norm are approximately $2$ and $4$ for the linear and quadratic Lagrange elements, respectively. These same convergence rates were obtained with a manufactured solution in section \ref{SubsectionManufSol} (see Figure \ref{Fig:L2errorH1errorP1P2elementsDirichlet}). Formula \eqref{ConvergenceRateFormula}, used on Table \ref{tab:ExactSolL2H1errorP1P2Neumann}, gives approximately the same rates of convergence as in Figure \ref{Fig:L2errorH1errorP1P2elementsDirichletVsNeumann} (for Neumann boundary conditions). While the quadratic Lagrange elements provide smaller error than the linear Lagrange elements, they are also the most computationally expensive. However, the quadratic Lagrange elements with mesh size $2h$ provide superior accuracy than the linear elements with mesh size $h$, and both require approximately the same computational cost. Therefore, high-order finite element approximations are beneficial in terms of efficiency. We also have computed the $L^2$- and $H^1$-error on the auxiliary $w_h$ and have found the same orders of convergence as for $u_h$, except in the $H^1$-error for quadratic Lagrange elements, which leads to super-convergence of order $\mathcal{O}(h^3)$.
\begin{table}[h!]
\centering
\begin{footnotesize}
\begin{tabular}{l|c|c|c|c|c|c|l}
\cline{2-7}
&\multicolumn{3}{c|}{Linear Lagrange elements} & \multicolumn{3}{c|}{Quadratic Lagrange elements} \\ \cline{1-7}
\multicolumn{1}{ |c| }{Mesh} & $L^2$-error& $H^1$-error & CPU (s) & $L^2$-error& $H^1$-error & CPU (s)  \\ \cline{1-7}
\multicolumn{1}{ |c| }{$25\times 25$}&$14.4303\times 10^{-6}$&$0.223426$&$15.80$&$599.199\times 10^{-9}$&$41.7231\times 10^{-4}$ &$25.98$   \\ \cline{1-7}
\multicolumn{1}{ |c| }{$50\times 50$}&$3.64952\times 10^{-6}$&$0.111751$&$25.99$&$37.8744\times 10^{-9}$&$10.4412\times 10^{-4}$  &$63.63$   \\ \cline{1-7}
\multicolumn{1}{ |c| }{$100\times 100$}&$0.91474\times 10^{-6}$&$0.055880$&$63.36$& $2.37924\times 10^{-9}$&$2.61096\times 10^{-4}$&$221.69$  \\ \cline{1-7}
\multicolumn{1}{ |c| }{$200\times 200$}&$0.22881\times 10^{-6}$&$0.027940$&$222.52$& $0.14958\times 10^{-9}$&$0.65278\times 10^{-4}$&$971.67$ \\ \cline{1-7}
\end{tabular}
\caption{\footnotesize CPU time, $L^2$- and $H^1$-error as  a function of the element size $h$, computed with linear and quadratic Lagrange elements, respectively. Results are obtained with Dirichlet boundary conditions.}
\label{tab:ExactSolL2H1errorP1P2Dirichlet}
\end{footnotesize}
\end{table}
\begin{table}[h!]
\centering
\begin{footnotesize}
\begin{tabular}{l|c|c|c|c|c|c|l}
\cline{2-7}
&\multicolumn{3}{c|}{Linear Lagrange elements} & \multicolumn{3}{c|}{Quadratic Lagrange elements} \\ \cline{1-7}
\multicolumn{1}{ |c| }{Mesh} & $L^2$-error& $H^1$-error & CPU (s) & $L^2$-error& $H^1$-error & CPU (s)  \\ \cline{1-7}
\multicolumn{1}{ |c| }{$25\times 25$}&$204.486\times 10^{-6}$&$0.223103$&$16.10$&$3799.57\times 10^{-9}$&$41.0548\times 10^{-4}$ &$27.01$   \\ \cline{1-7}
\multicolumn{1}{ |c| }{$50\times 50$}&$51.5883\times 10^{-6}$&$0.111696$&$26.21$&$336.114\times 10^{-9}$&$10.3552\times 10^{-4}$  &$64.37$   \\ \cline{1-7}
\multicolumn{1}{ |c| }{$100\times 100$}&$12.9592\times 10^{-6}$&$0.055871$&$63.97$& $29.7100\times 10^{-9}$&$2.60006\times 10^{-4}$&$223.69$  \\ \cline{1-7}
\multicolumn{1}{ |c| }{$200\times 200$}&$3.24528\times 10^{-6}$&$0.027939$&$224.72$& $2.62640\times 10^{-9}$&$0.65140\times 10^{-4}$&$992.07$ \\ \cline{1-7}
\end{tabular}
\caption{\footnotesize CPU time, $L^2$- and $H^1$-error as a function of the element size $h$, computed with linear and quadratic Lagrange elements, respectively. Results are obtained with Neumann boundary conditions.}
\label{tab:ExactSolL2H1errorP1P2Neumann}
\end{footnotesize}
\end{table}

For the error and convergence analysis in time, we compute a reference solution ($u_h^{\text{ref}}$,$w_h^{\text{ref}}$) on  a grid of $100\times 100$ using linear Lagrange elements and a time step $\Delta t= 10^{-5}$. We then analyze the numerical convergence in time by varying the time step $\Delta t$ systematically between $0.001$ and $0.008$. For each time step, the errors on $u_h$ and $w_h$, computed as 
\begin{equation}
\dfrac{\Vert u_h^{\text{ref}} -u_h \Vert_{L^2(\Omega)}}{\Vert u_h^{\text{ref}}\Vert_{L^2(\Omega)}}\quad \text{ and }\quad \dfrac{\Vert w_h^{\text{ref}} -w_h \Vert_{L^2(\Omega)}}{\Vert u_h^{\text{ref}}\Vert_{L^2(\Omega)}},
\end{equation}
respectively, are recorded at the final time $T=0.201$. Table~\ref{L2errorExactSolTimeOrderSBDF2} shows the error, computational time and order of convergence of the numerical solution with respect to the reference solution. The rate of convergence in time is calculated using the formula \eqref{ConvergenceRateFormula}, where $\tau$ is the refinement factor between consecutive time steps. Here, also $\tau=2$ was chosen in our numerical tests. A second-order convergence rate  in time is obtained.
\begin{table}[h!]
\centering
\begin{tabular}{c|c|c|c|c|c|c}
\cline{2-5}
\multicolumn{1}{c|}{}&\multicolumn{2}{|c|}{Relative $L^2$-error}&\multicolumn{2}{|c|}{Order of convergence}&\multicolumn{1}{|c}{} \\ \cline{1-6}
\multicolumn{1}{|c|}{$\Delta t$}& Error on $u_h$&Error on $w_h$ & Order for $u_h$ & Order for $w_h$ & CPU (s) \\ \cline{1-6}
\multicolumn{1}{|c|}{$0.008$}&$1.4340\times 10^{-6}$&$3.4856\times 10^{-5}$&$2.079$&$2.079$&$50.367$ \\ \cline{1-6}
\multicolumn{1}{|c|}{$0.004$}& $3.3927\times 10^{-7}$&$8.2465\times 10^{-6}$&$2.038$&$2.038$&$95.467$ \\ \cline{1-6}
\multicolumn{1}{|c|}{$0.002$}& $8.2570\times 10^{-8}$& $2.0070\times 10^{-6}$&$2.019$&$2.019$&$190.636$ \\ \cline{1-6}
\multicolumn{1}{|c|}{$0.001$}& $2.0369\times 10^{-8}$& $4.9512\times 10^{-7}$& --- & ---  &  $386.169$ \\ \cline{1-6}
\end{tabular}
\caption{Relative $L^2$-error as a function of time step $\Delta t$, order of convergence (in time) and CPU time. Results are obtained with Dirichlet boundary conditions.}
\label{L2errorExactSolTimeOrderSBDF2}
\end{table}

\subsection{Numerical results --- Applications}
\subsubsection{Lubrication}
\begin{figure}[!h]
\centering
\begin{tabular}{cccc}
\includegraphics[scale=0.2]{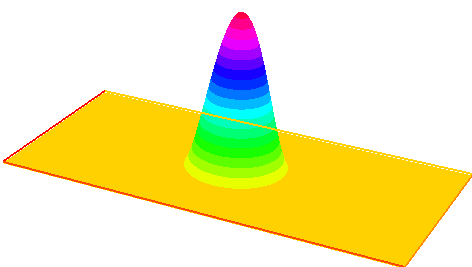}&\includegraphics[scale=0.2]{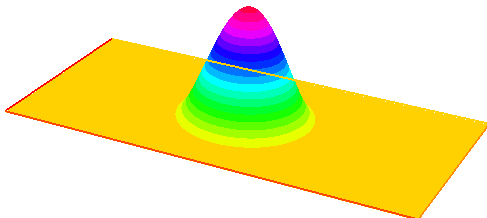}&\includegraphics[scale=0.2]{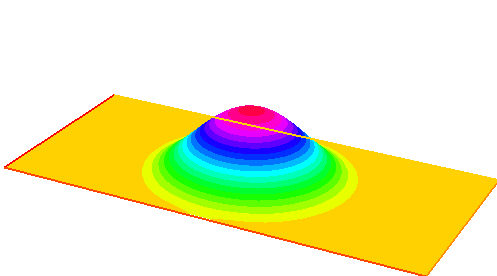}&\includegraphics[scale=0.2]{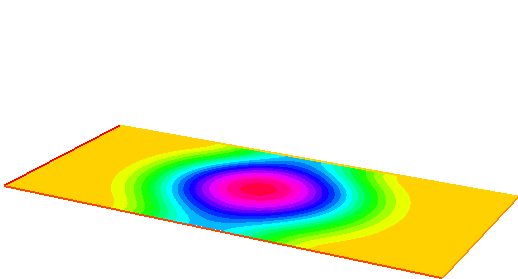}\\
\includegraphics[scale=0.3]{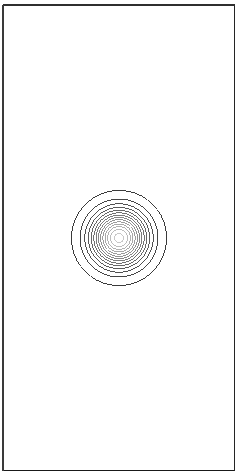}&\includegraphics[scale=0.3]{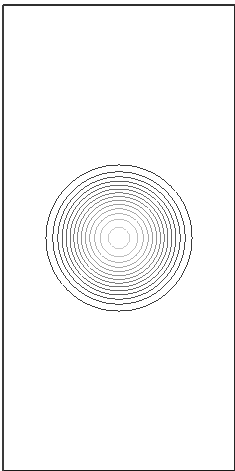}&\includegraphics[scale=0.3]{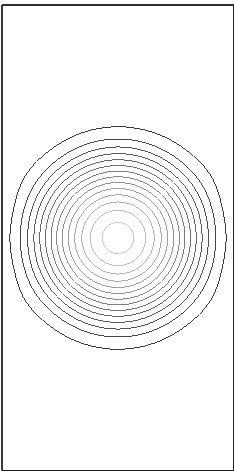}&\includegraphics[scale=0.3]{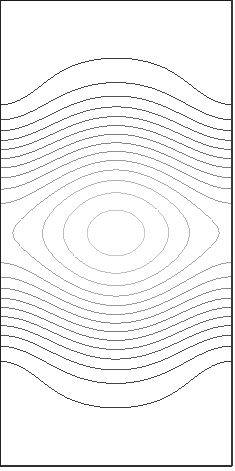}\\
$t=0$&$t=0.0001$&$t=0.001$&$t=0.01$
\end{tabular}
\caption{Numerical solution of \eqref{lubSyst} starting from \eqref{LubInitCond}. Top: A 3D view of the droplet spreading on a solid surface. Bottom: Contours at various times of the simulated solution. $\Delta t=10^{-5}$.}
\label{Fig:Contours3dViewDropletSpreading}
\end{figure}

Equation \eqref{lubSyst} appears in the lubrication approximation of a thin viscous film, which is driven by surface tension. In this context, $u$ describes the height of the liquid film which spreads on a solid surface \cite{greenspan1978,Hocking1981}. The numerical solutions of high-order nonlinear degenerate diffusion equations like  \eqref{lubSyst} require special numerical techniques in order to preserve the non-negativity of solutions.  Constructing finite element methods that guarantee the non-negativity of the solutions of \eqref{lubSyst}, has been initiated by Barrett et al.\  \cite{Barrett1} by enforcing at each time step the non-negativity of the approximate solution as a constraint, which leads to a variational inequality. The solutions of the one-dimensional version of \eqref{lubSyst} using a non-negativity preserving scheme was carried out in \cite{Grun2000}. The solution of a modified problem, which consists in approximating the mobility term $f(u)=u$ by
\begin{equation}
f_{\xi}(u)=\dfrac{u^5}{\xi u + u^4},
\label{RegularizationTerm}
\end{equation}
has been shown by \cite{BERNIS1990} to remain positive. This regularization is used in \cite{Bertozzi1,WITELSKI2003} to compute numerical solutions of \eqref{lubSyst} for the initial data
\begin{equation}
u(x,y,0)=\delta+Ce^{-\sigma(x^2+y^2)},
\label{LubInitCond}
\end{equation}
where $\delta>0$ represents the thickness of an ultra-thin precursor film under the Gaussian fluid droplet centered at the origin.   
Here, we solve \eqref{lubSyst} without regularization, i.e.\ with the degenerate mobility function $f(u)=u$. The computational domain is the rectangle $[-0.5,0.5]\times[-1,1]$. Homogeneous Neumann conditions \eqref{HomoNeumannBoundCond} are imposed on the boundary.  We use linear Lagrange elements on a uniform mesh $70\times 140$, i.e.\ a mesh size $h=0.02$. We set $\delta=0.01$, $\sigma=80$ and $C=2$. The same values for these parameters and mesh size were used in \cite{Bertozzi1}. The contours of the numerical solutions and a 3D view of the droplet  at different times are shown in Figure~\ref{Fig:Contours3dViewDropletSpreading}. 
The algorithm \eqref{SBDF2ThinWeakForm}-\eqref{SBDF1ThinWeakForm} preserves the  positivity of the solution up to a maximal time step of $\Delta t\simeq 10^{-5}$. We have also performed numerical simulations for $C=1$ as in \cite{WITELSKI2003}, for which the positivity has been preserved up to a maximal time step of $\Delta t\simeq 0.01$, while $\Delta t=10^{-5}$  was found  in \cite{WITELSKI2003} as the maximal time step to preserve the positivity.

\begin{figure}[!h]
\centering
\begin{tabular}{cc}
\includegraphics[scale=0.57]{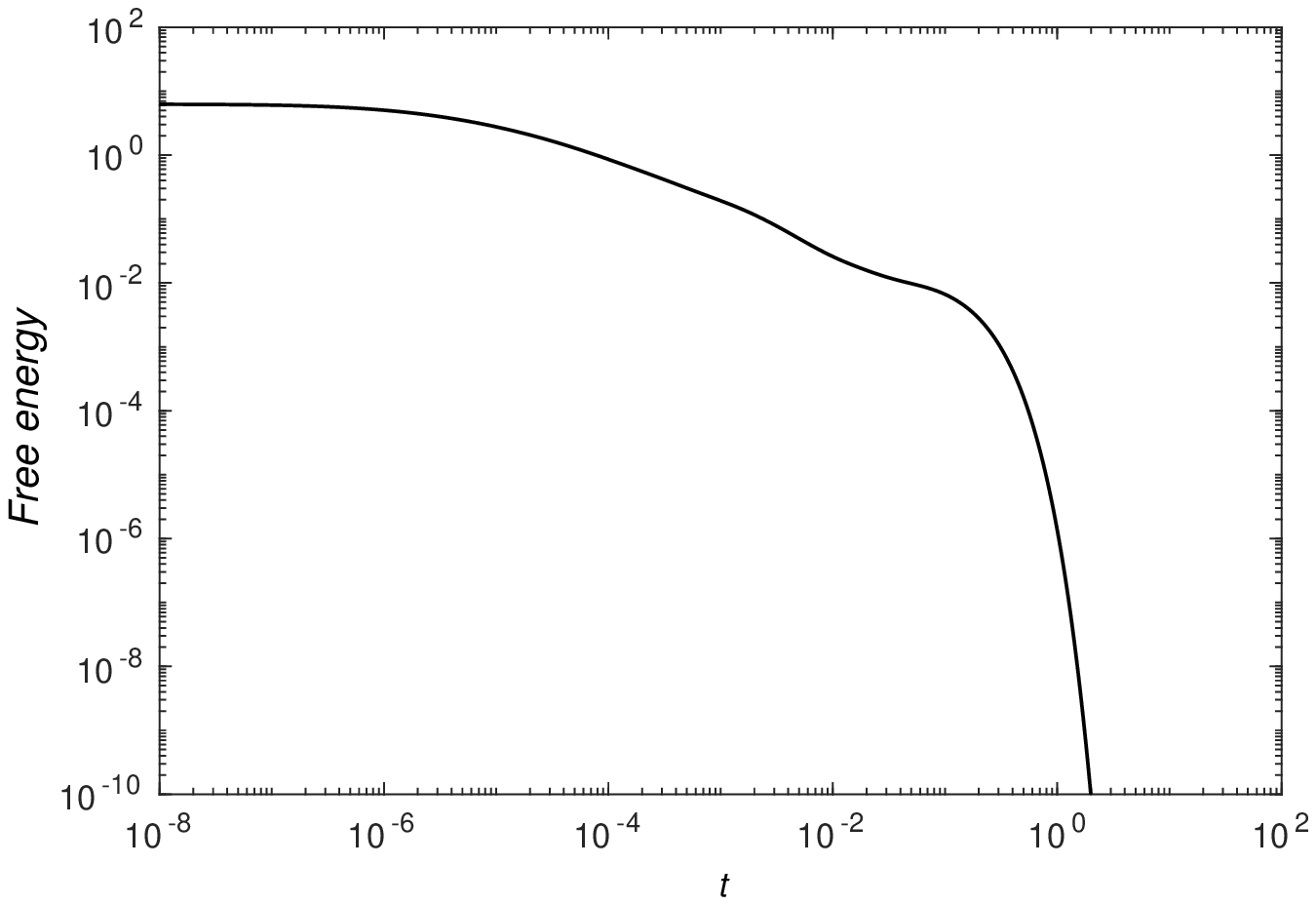}&\includegraphics[scale=0.57]{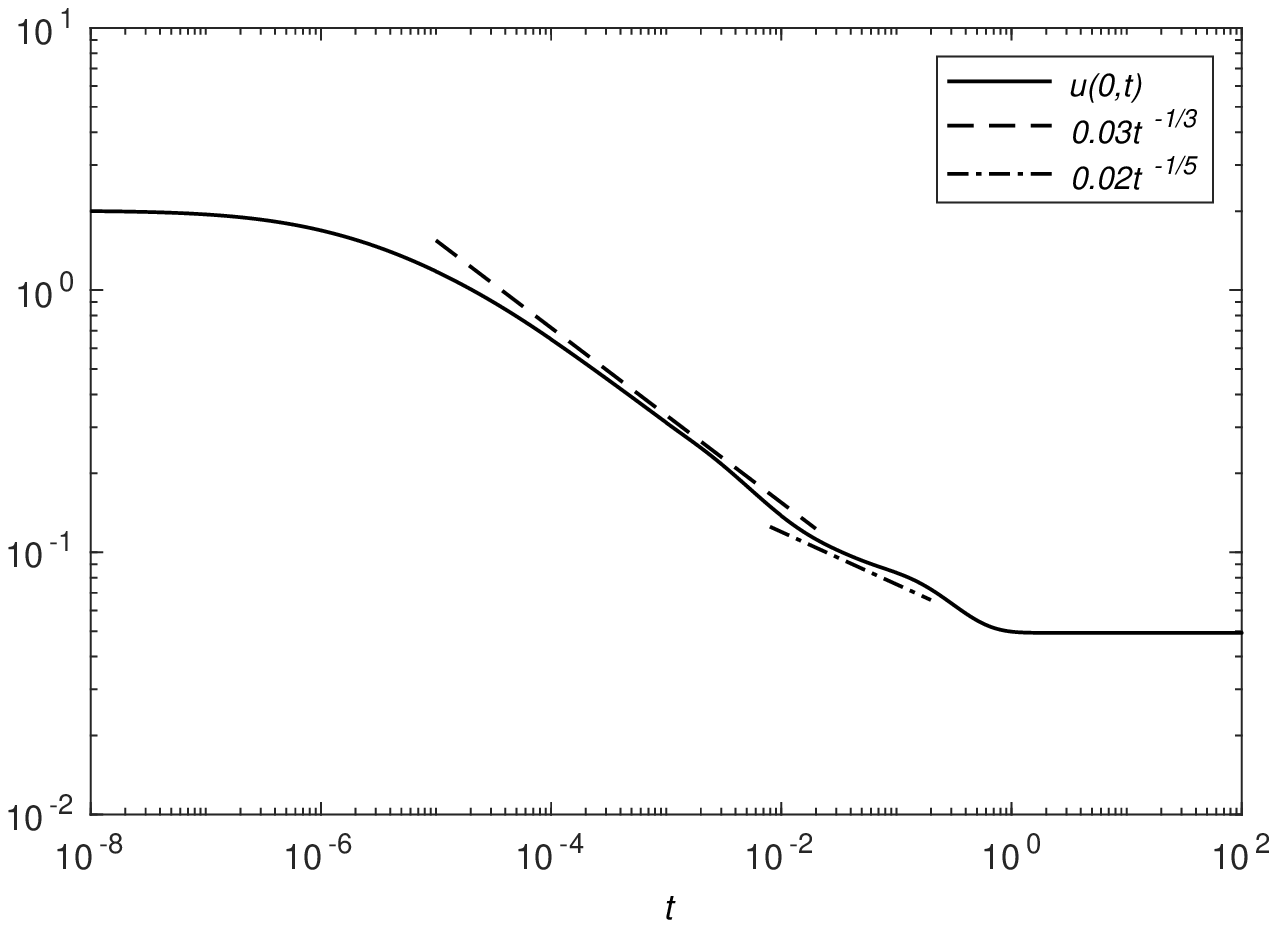}\\
$(a)$ & $(b)$
\end{tabular}
\caption{$(a)$ Monotonic decrease of the total free energy estimate. $(b)$ Evolution of the droplet maximal height $u(0,t)$ showing the intermediate asymptotics.}
\label{Fig:FreeEnergyandHeight}
\end{figure}

To show that the numerical scheme satisfies the energy decreasing property \eqref{EnergyProperty}, we define the discrete energy functional of \eqref{lubSyst} as
\begin{equation}
\varepsilon_h(u_h^n)=\frac{1}{2}\int_{\Omega}\vert\nabla u_h^n\vert^2\,d\bm{x},
\end{equation}
and consider the evolution of the discrete energy. Figure~\ref{Fig:FreeEnergyandHeight}$(a$) shows the time-evolution of the discrete total energy. Our results demonstrate that the total discrete energy is non-increasing which agrees well with the total energy dissipation property. 

The maximum of the solution of \eqref{lubSyst} with the initial condition \eqref{LubInitCond} occurs at the origin for all times. Figure~\ref{Fig:FreeEnergyandHeight}$(b)$ shows the height of the droplet at the origin as time evolves. After an initial transient phase ($t>10^{-5})$, the maximal height $u(0,t)$ evolves in $t^{-1/3}$ for some times, and in $t^{-1/5}$ for longer times, before approaching a uniform flat state for large times ($t>10^{-1}$). The same evolutions are observed in \cite{Bertozzi1,WITELSKI2003}. 

\subsubsection{Phase separation}
Spinodal decomposition is a mechanism by which a solution of two or more components can separate into distinct phases. Spinodal decomposition in a binary mixture is well-described by the Cahn-Hilliard equation \cite{CAHN1962,CahnHilliard1958}. Here,  we are interested in the phase separation of a homogeneous binary mixture governed by \eqref{ThinFilmSyst} to further test our method \eqref{SBDF2ThinWeakForm}-\eqref{SBDF1ThinWeakForm}. As in \cite{Vignal2017}, we choose a mobility function $f(u)=(1-u^2)$, a free energy function $\varphi(u)=\frac{1}{4\epsilon^2}(u^2-1)^2$ and set the parameters $\gamma=1$ and $\epsilon=0.03$.  In this case, equation \eqref{ThinFilmSyst} describes the state of a relative concentration $u=(2\frac{u_a}{u_a+u_b}-1)$, where $u_a$ and $u_b$ are the concentrations of the two components \cite{GARCKE2003}. With this double well form of the energy function, the system is locally in one of the two phases if $u$ is close to one of the two minima $\pm 1$ of $\varphi$ \cite{Elliott89}.  

\begin{figure}[!h]
\centering
\begin{tabular}{cc}
\includegraphics[scale=0.2]{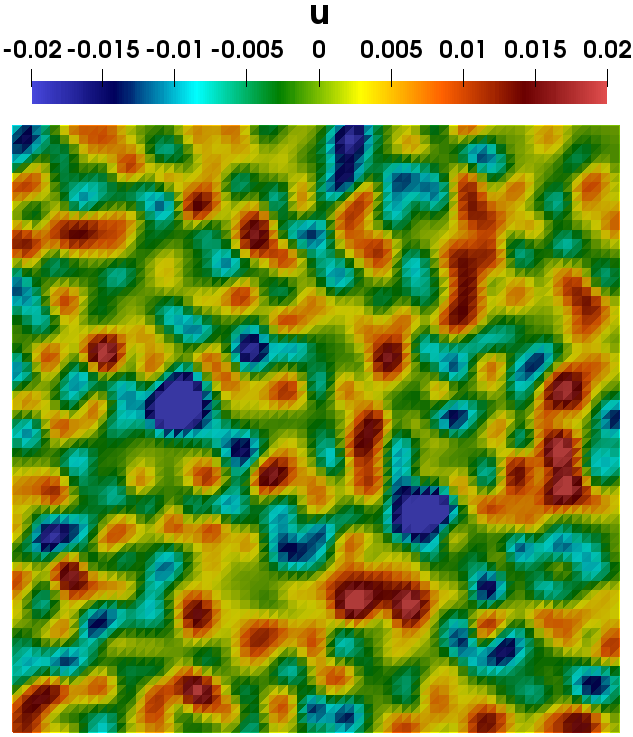}&\includegraphics[scale=0.2]{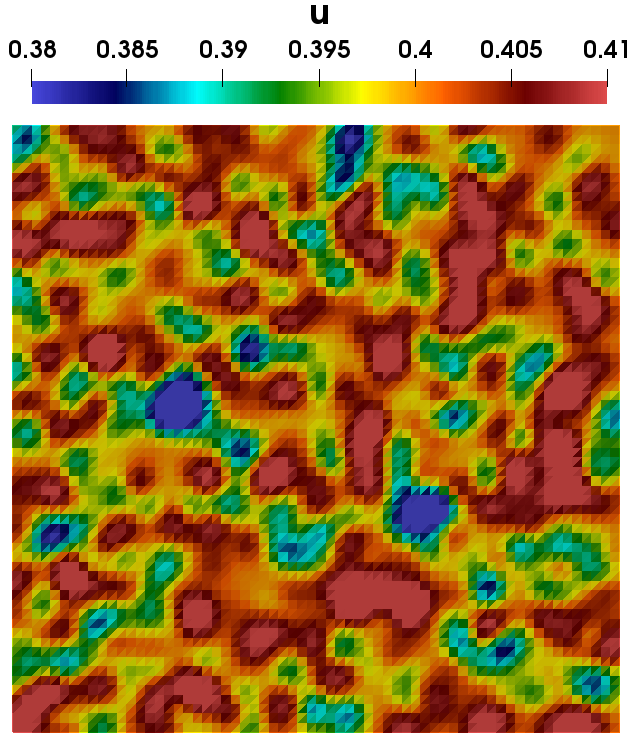}\\
$\overline{u}=0$&$\overline{u}=0.4$
\end{tabular}
\caption{Initial solution for different values of $\overline{u}$.}
\label{Fig:PhaseSeparationInitSol}
\end{figure}

The initial condition is chosen as
\begin{equation}
u_0(\bm{x})=\overline{u}+ \sum_{i=1}^{M} c_ir_i(\bm{x},\bm{a}_i),
\label{randInitData} 
\end{equation}
where $\overline{u}\in[-1,1]$, $\bm{a}_i\in\mathbb{R}^2$, $c_i\in\mathbb{R}$, for $i=1,...,M$, and $r_i(\bm{x},\bm{a}_i)=e^{-1000\vert\bm{x}-\bm{a}_i\vert^2}$. We have randomly generated $M=1000$ points $\bm{a}_i$ in the square domain $[0,1]\times[0,1]$, and $1000$ scalars $c_i\in[-0.01,0.01]$. All numerical results in this section are computed on a $64\times 64$ grid with quadratic Lagrange elements. Figure~\ref{Fig:PhaseSeparationInitSol} shows the initial condition for two different values of $\overline{u}$. 

\begin{figure}[!h]
\centering
\begin{tabular}{cccc}
\includegraphics[scale=0.15]{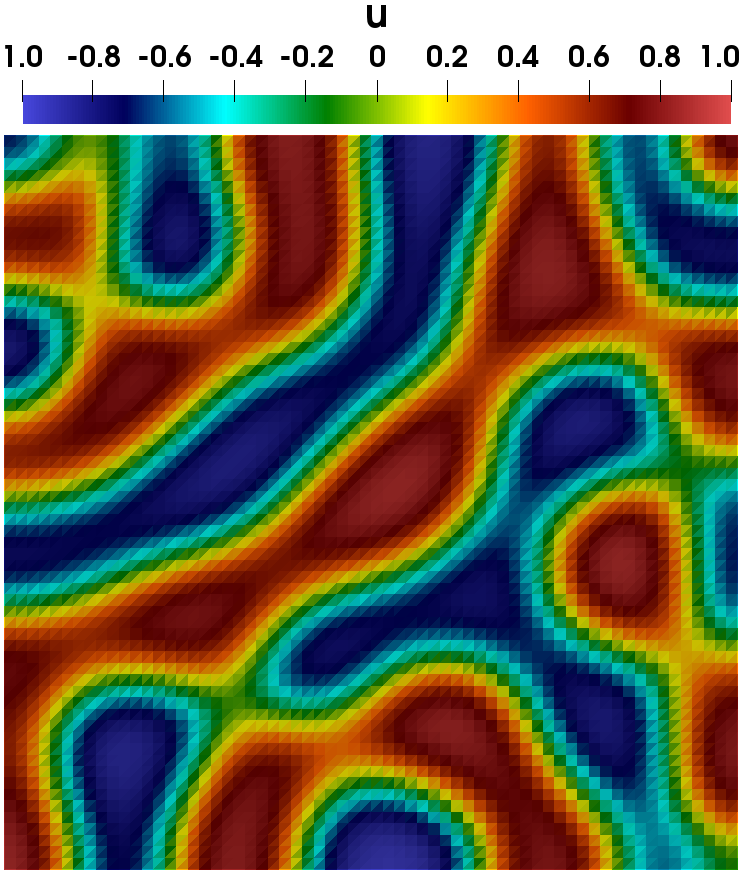}&\includegraphics[scale=0.15]{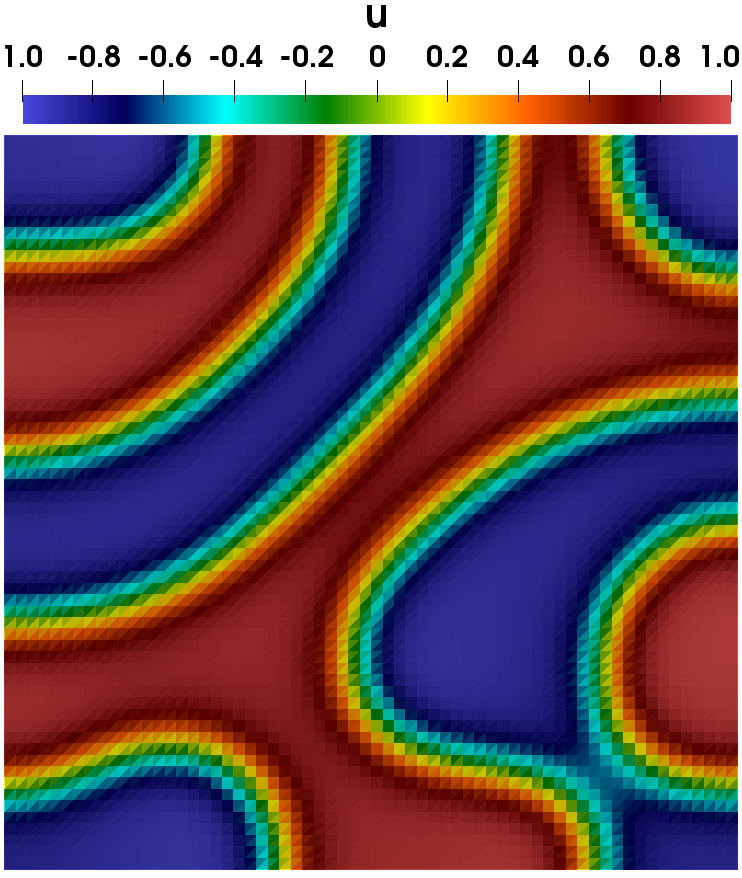}&\includegraphics[scale=0.15]{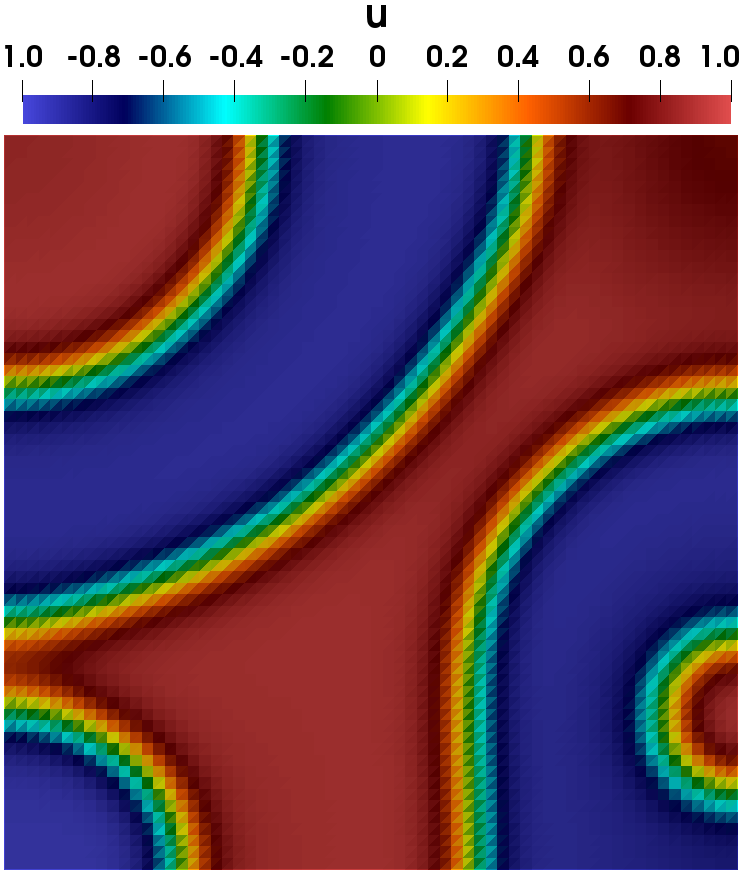}&\includegraphics[scale=0.15]{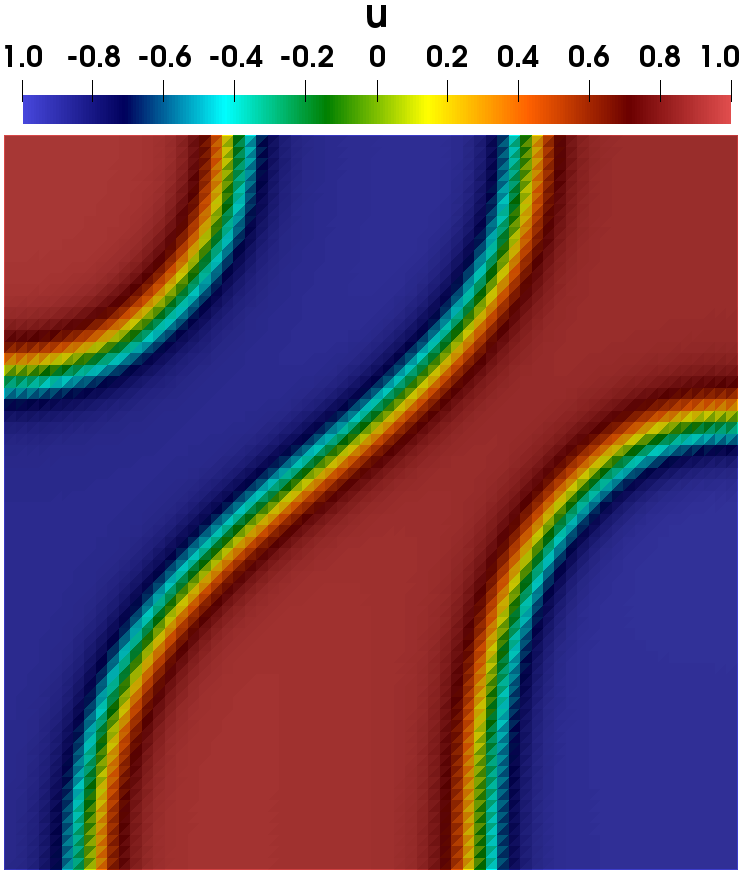}\\
$t=5\times 10^{-5}$&$2\times 10^{-4}$&$t=1\times 10^{-3}$&$t=2\times 10^{-3}$
\end{tabular}
\caption{Simulation of a spinodal decomposition starting from a uniformly distributed random initial condition. $\Delta t=10^{-6}$.}
\label{Fig:SpinodalDecomp}
\end{figure}

We first consider the initial solution for $\overline{u}=0$. The dynamics of the binary fluid are shown in Figure~\ref{Fig:SpinodalDecomp}. From its initial concentration, the system separates quickly into two phases whose composition is determined by the two minima $\pm 1 $ of the free energy function. Small regions of the same phase combine together to form larger patterns, reducing the interface boundaries between the two phases.

\begin{figure}[!h]
\centering
\begin{tabular}{cccc}
\includegraphics[scale=0.15]{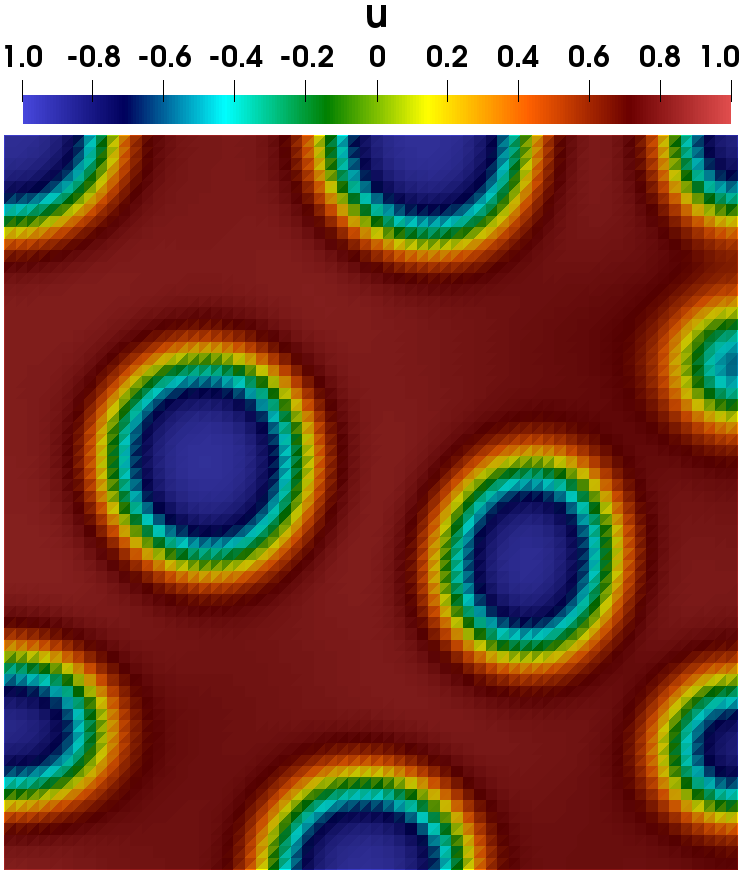}&\includegraphics[scale=0.15]{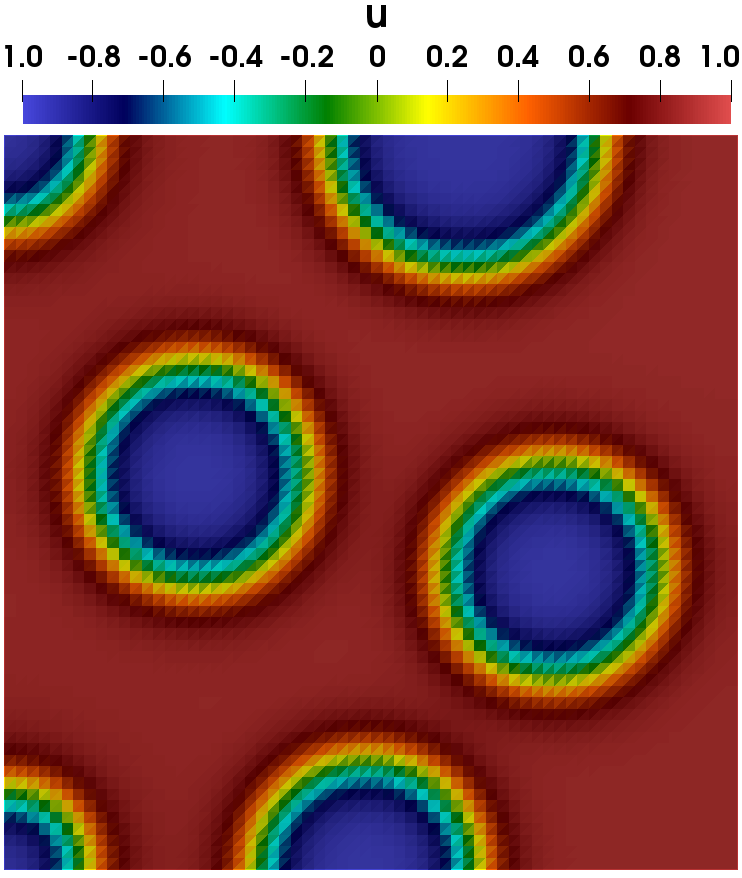}&\includegraphics[scale=0.15]{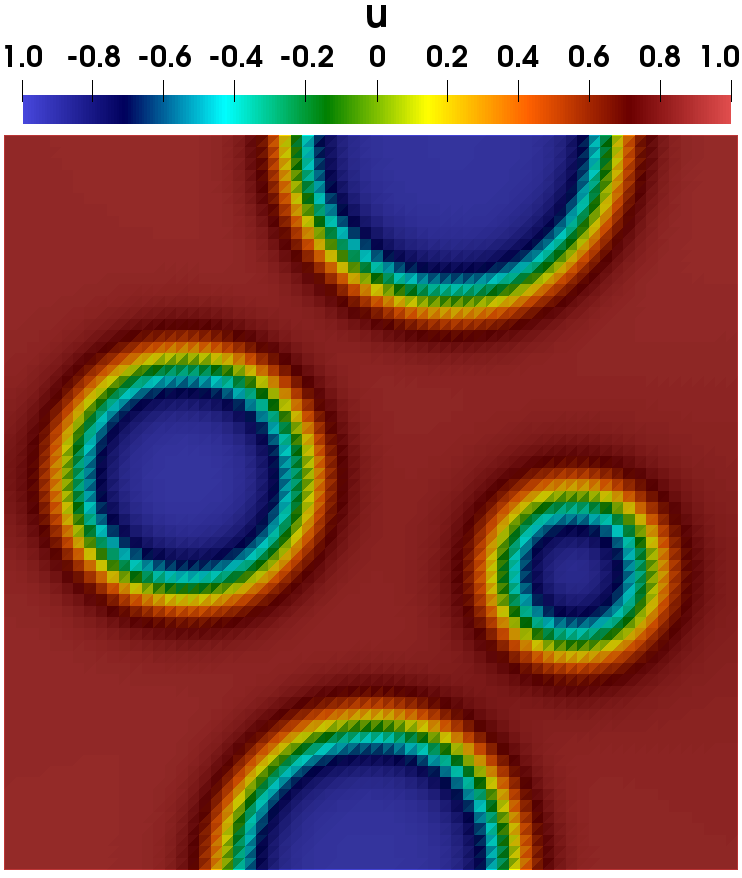}&\includegraphics[scale=0.15]{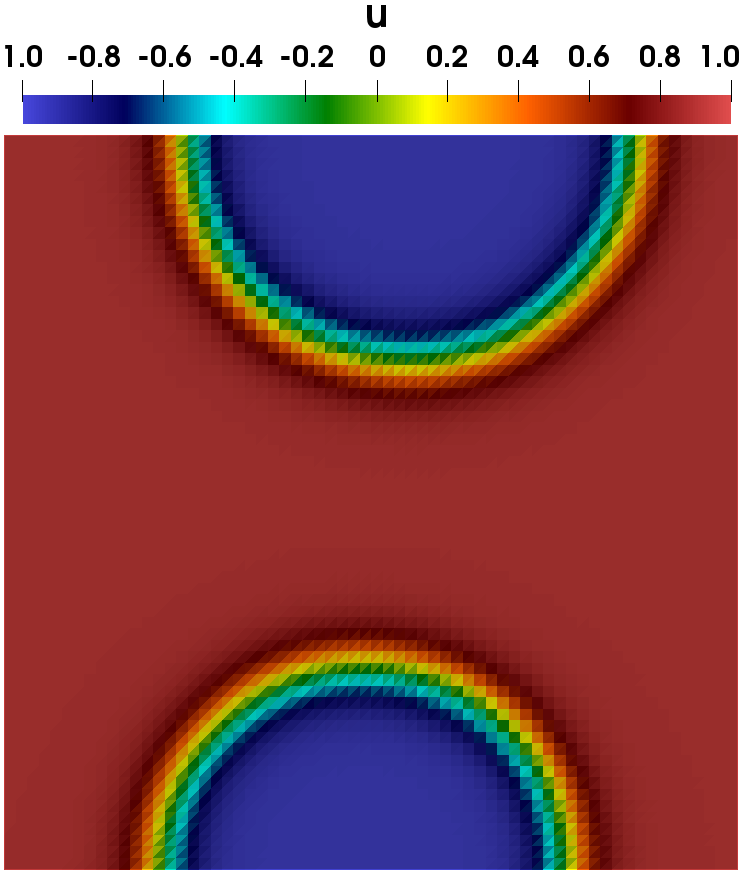}\\
$t=1.25\times 10^{-4}$&$1\times 10^{-3}$&$t=2.5\times 10^{-3}$&$t=1\times 10^{-2}$
\end{tabular}
\caption{Simulation of a Ostwald ripening process  starting from a perturbation of an uniform state. $\Delta t=10^{-6}$.}
\label{Fig:OstwaldRipening}
\end{figure}

A particular regime of interest is when one phase is significantly more abundant in the initial mixture. In such case, the minority phase emerges and forms droplets. Larger particles grow at the expense of the smaller ones which
shrink and disappear. This form of competitive growth is known as the Ostwald ripening \cite{Pismen2006}. To illustrate such regime, we use the initial solution for $\overline{u}=0.4$. Figure~\ref{Fig:OstwaldRipening} shows the time-evolution of the binary mixture computed. We can observe the formation of droplets. Smaller droplets are absorbed into larger ones through diffusion, exhibiting clearly the phenomenon of Ostwald ripening. 

\begin{figure}[!h]
\centering
\begin{tabular}{cc}
\includegraphics[scale=0.58]{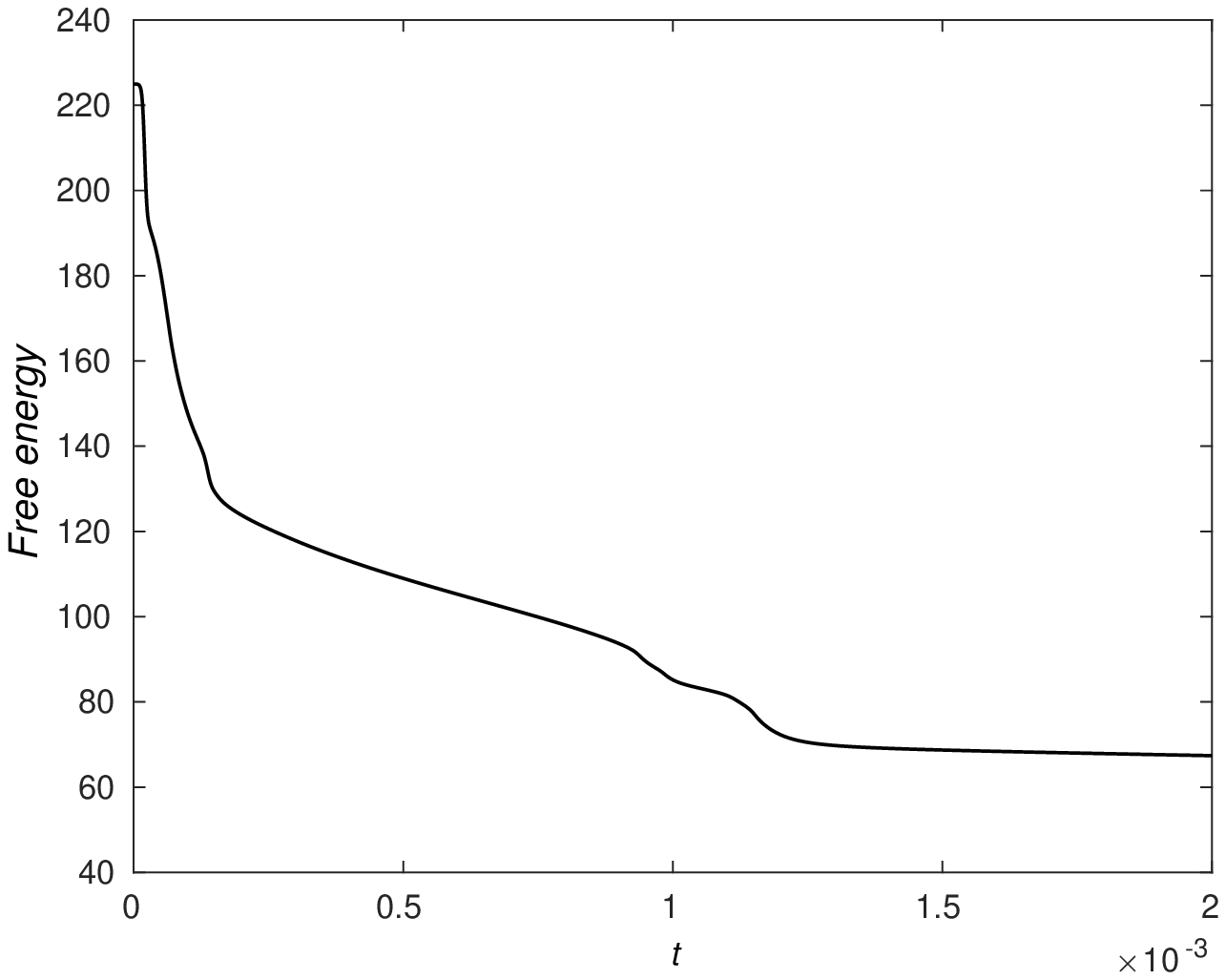}&\includegraphics[scale=0.58]{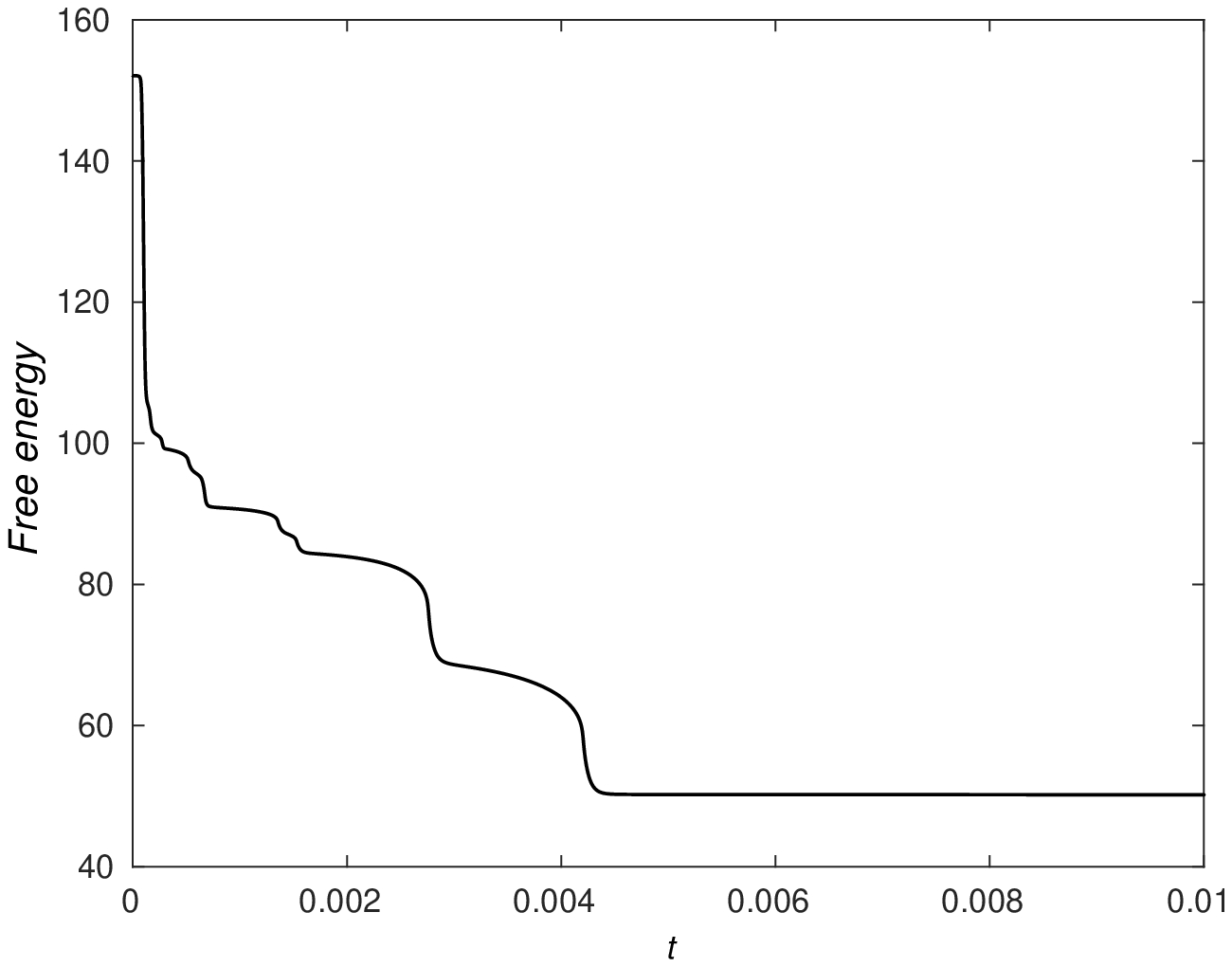}\\
$(a)$&$(b)$
\end{tabular}
\caption{Free energy evolutions. $(a)$ for the spinodal decomposition  and $(b)$ for the Ostwald ripening.}
\label{Fig:EnergyPlot}
\end{figure}

\begin{figure}[!h]
\centering
\begin{tabular}{cc}
\includegraphics[scale=0.58]{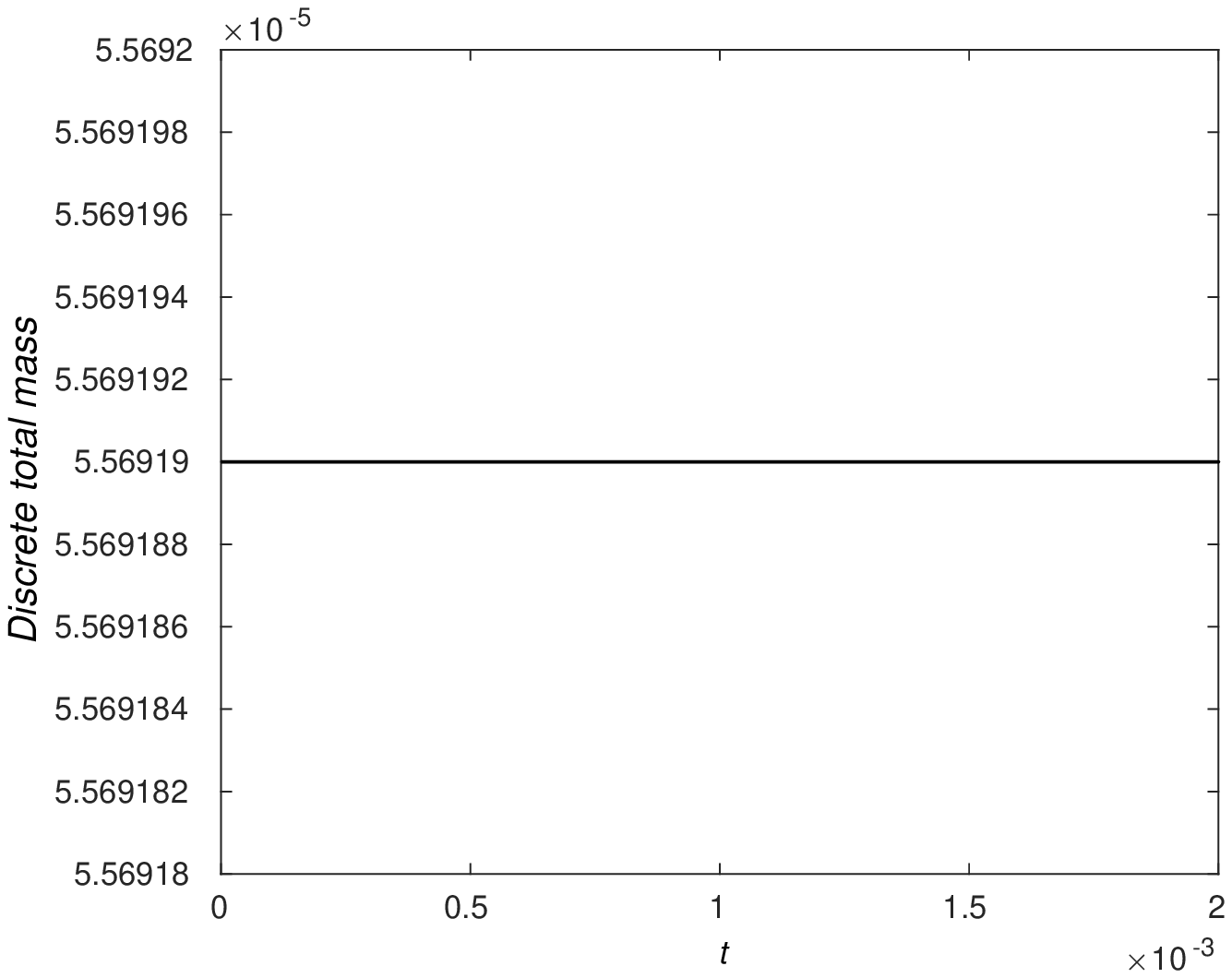}&\includegraphics[scale=0.58]{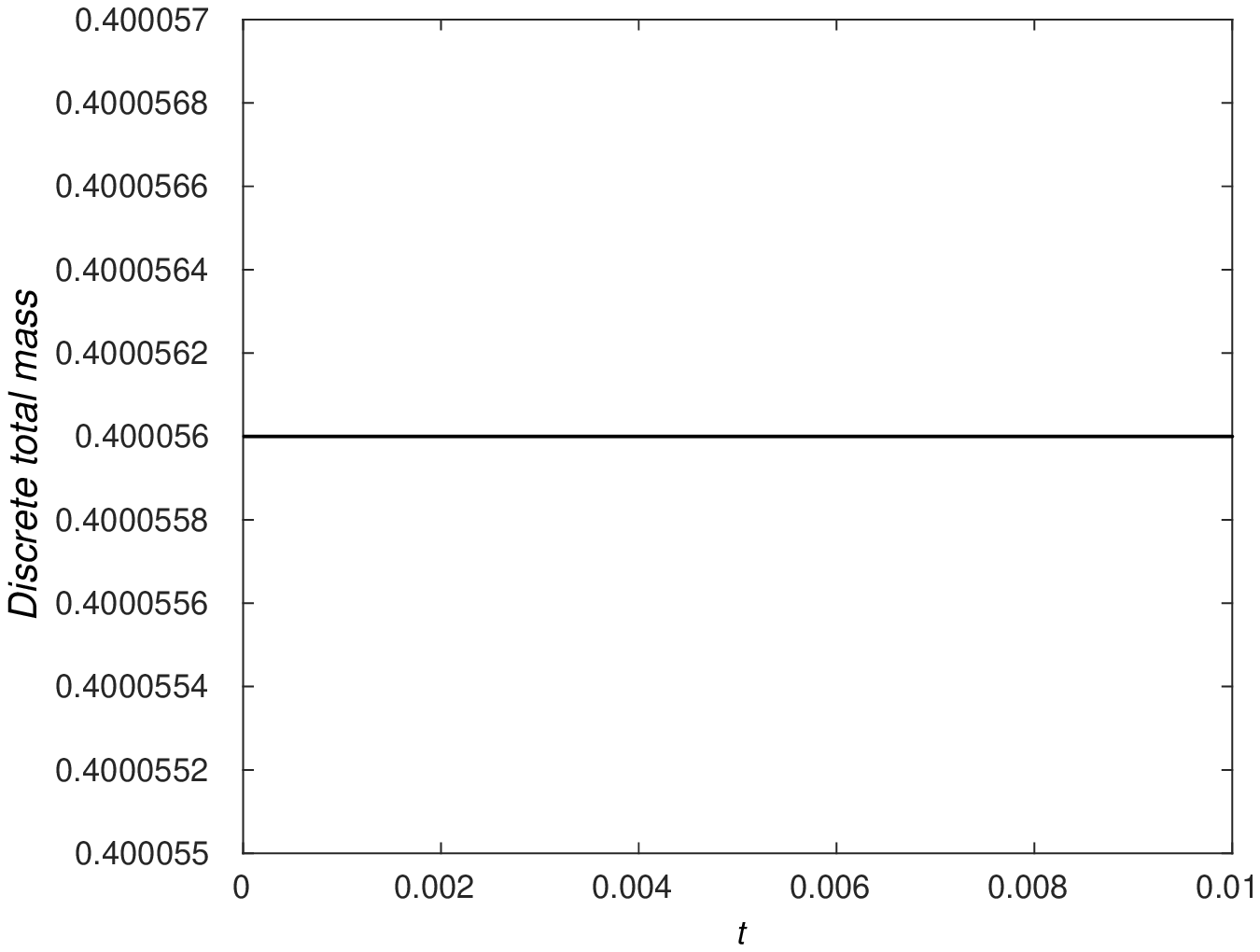}\\
$(a)$&$(b)$
\end{tabular}
\caption{Numerical verification of the total mass conservation over time. $(a)$ for the spinodal decomposition  and $(b)$ for the Ostwald ripening.}
\label{Fig:MassPlot}
\end{figure}

We qualitatively observe the same behavior of the solution as in \cite{Vignal2017}. The plots of the energy evolution in Figure~\ref{Fig:EnergyPlot} illustrates that the total discrete energy satisfies the energy dissipation property. Plots of the discrete total mass are shown in Figure~\ref{Fig:MassPlot}, which  illustrate the mass-conservative property of the linear scheme.
 We observe that the magnitude and the nonlinearity of the mobility have a significant direct impact on the stability of the scheme. We have noticed from numerical experiments that a small time step is required for the stability of numerical solutions involving non constant mobility functions. The time step could be taken larger for constant mobility functions and regular energy functions. The parameter $\gamma$ has an impact on the length scale of the diffuse interface by balancing the surface energy term with the bulk free energy \cite{Bertozzi4}. In terms of the stability of the scheme, there is a strong link between the parameter $\gamma$ and the free energy term. 

\subsubsection{Electrowetting}
Electrowetting has become one of the most widely used tools for manipulating tiny amounts of liquids on surfaces. It consists in the modification of the wetting properties of a surface with an applied electric field. A diffuse-interface model for drop motion, due to electrowetting, in a Hele-Shaw geometry is proposed in \cite{Bertozzi2}. The model is governed by the Cahn-Hilliard model with a spatially dependent energy
\begin{equation}
\left\lbrace
\begin{aligned}
&\epsilon u_t-\nabla\cdot(u\nabla w)=0,\\
&w=-\epsilon^2\Delta u +\varphi'(u) -\epsilon\rho(\bm{x}),
\end{aligned}
\right.
\label{WettingSyst}
\end{equation}
with a double well energy function $\varphi(u)=(u-\epsilon)^2(u-1)^2$, where $\rho(\bm{x})=\lambda\chi(\bm{x})$ defines a local energy (electric field) with characteristic function $\chi(\bm{x})$ and strength $\lambda$. The local energy  term introduces the electrowetting into the diffuse-interface model. The parameter  $\epsilon=0.0427$ controls the diffuse interface thickness. For more details and the experimental setup, we refer to \cite{Bertozzi2}.

We solve this problem using a Gaussian initial data
\begin{equation}
u_0(\bm{x})=\delta+e^{-10\vert \bm{x}-\bm{a} \vert^2},
\end{equation}
representing a droplet centered at the point $\bm{a}$ that will be specified for each test case. All numerical results in this section are computed on a $47\times 94$ grid in the rectangular domain $[-0.5,0.5]\times[-1,1]$ using linear Lagrange elements. This gives rise to a mesh size $h\approx0.03$ as in \cite{Bertozzi1,Bertozzi2}. We use a precursor film $\delta=\epsilon$  and a fixed time step $\Delta t=0.001$ instead of adaptative time steps as in \cite{Bertozzi1}.

\begin{figure}[!h]
\centering
\begin{tabular}{ccccccc}
$(a)$\includegraphics[scale=0.2]{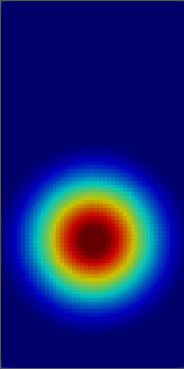}&\includegraphics[scale=0.2]{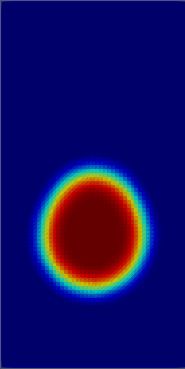}&\includegraphics[scale=0.2]{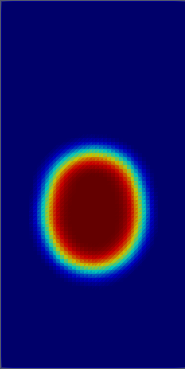}&\includegraphics[scale=0.2]{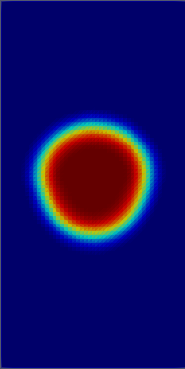}&\includegraphics[scale=0.2]{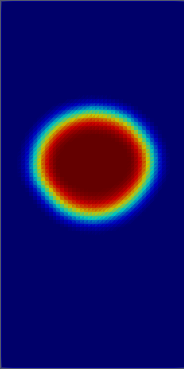}&\includegraphics[scale=0.2]{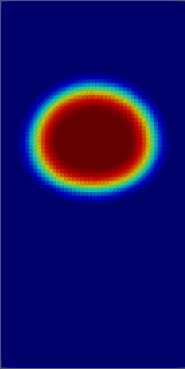}&\includegraphics[scale=0.2]{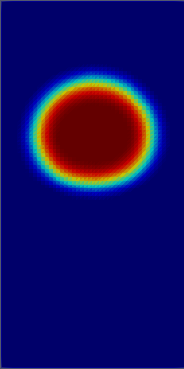}\\
$t=0$&$t=0.1$&$t=0.2$&$t=0.3$&$t=0.35$&$t=0.45$&$t=0.5$\\\\
$(b)$\includegraphics[scale=0.2]{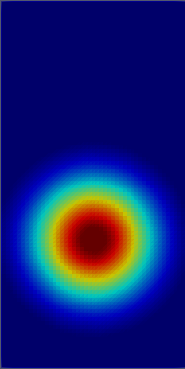}&\includegraphics[scale=0.2]{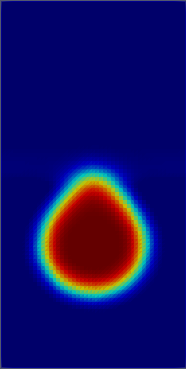}&\includegraphics[scale=0.2]{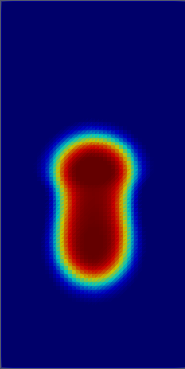}&\includegraphics[scale=0.2]{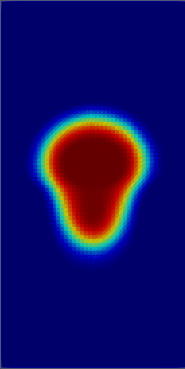}&\includegraphics[scale=0.2]{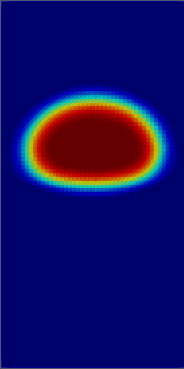}&\includegraphics[scale=0.2]{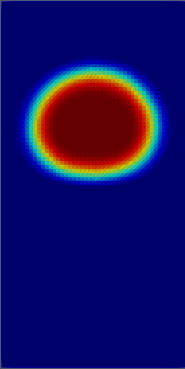}&\includegraphics[scale=0.2]{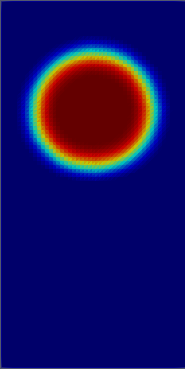}\\
$t=0$&$t=0.01$&$t=0.05$&$t=0.08$&$t=0.11$&$t=0.2$&$t=0.5$
\end{tabular}
\caption{Droplet translation by electrowetting: $(a)$ solution for $\lambda=0.75$ and $(b)$ $\lambda=3$. $\Delta t=10^{-3}$.}
\label{Fig:ElectrowettingDiffLambda}
\end{figure}

We first consider a droplet centered at point $\bm{a}=(0,-0.3)$ and then apply the electric field $\rho(\bm{x})$ on the surface $y>0$. The time-evolution of the droplet is displayed in Figure~\ref{Fig:ElectrowettingDiffLambda}. The surface of the droplet is deformed during its translation toward the region with low energy (with an electric field). The variation of the curvature of the droplet increases with the strength $\lambda$ of the electric field. The same evolutions are observed using numerical simulations in \cite{Bertozzi1,Bertozzi2} and experimentally in \cite{Bertozzi2}.

\begin{figure}[!h]
\centering
\begin{tabular}{ccccccc}
\includegraphics[scale=0.2]{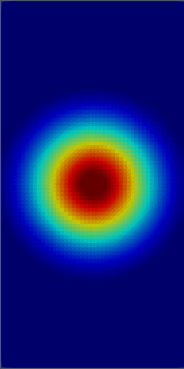}&\includegraphics[scale=0.2]{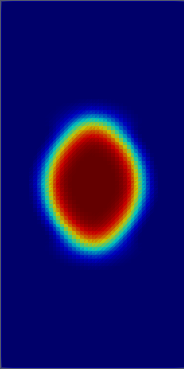}&\includegraphics[scale=0.2]{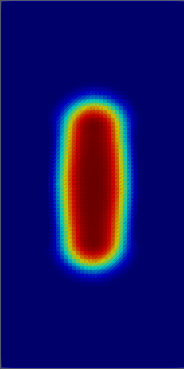}&\includegraphics[scale=0.2]{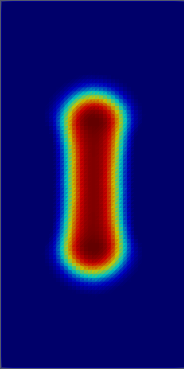}&\includegraphics[scale=0.2]{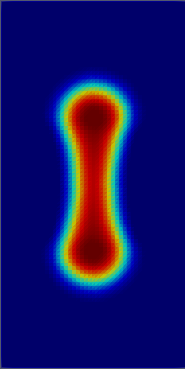}&\includegraphics[scale=0.2]{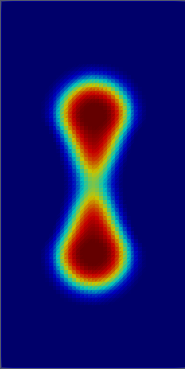}&\includegraphics[scale=0.2]{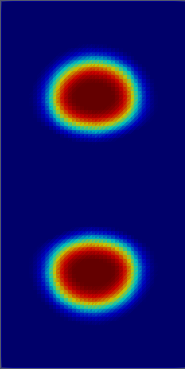}\\
$t=0$&$t=0.015$&$t=0.05$&$t=0.075$&$t=0.1$&$t=0.12$&$t=0.15$\\\\
\includegraphics[scale=0.2]{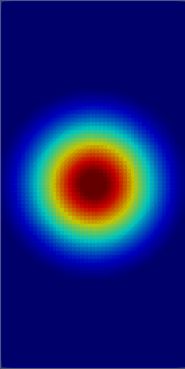}&\includegraphics[scale=0.2]{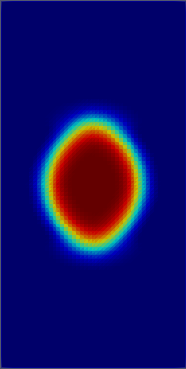}&\includegraphics[scale=0.2]{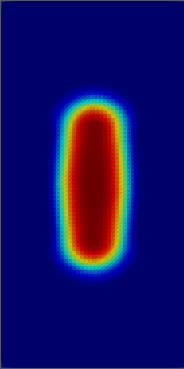}&\includegraphics[scale=0.2]{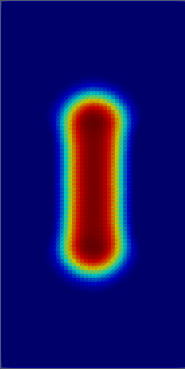}&\includegraphics[scale=0.2]{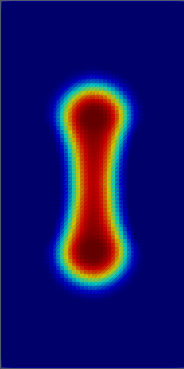}&\includegraphics[scale=0.2]{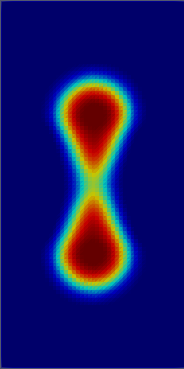}&\includegraphics[scale=0.2]{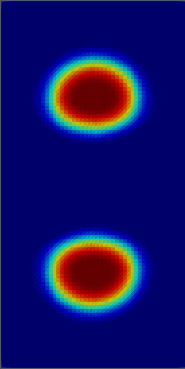}\\
$t=0$&$t=0.015$&$t=0.05$&$t=0.075$&$t=0.1$&$t=0.12$&$t=0.15$
\end{tabular}
\caption{Droplet splitting by electrowetting: (top) solution with the algorithm \eqref{SBDF2ThinWeakForm}-\eqref{SBDF1ThinWeakForm} and (bottom) solution with BDF2 scheme. $\Delta t=10^{-3}$, $\lambda=2$.}
\label{Fig3BDF2VsSBDF2}
\end{figure}

We next consider a droplet centered at the origin $\bm{a}=(0,0)$ and apply two electric fields of same strength $\lambda=2$ on the surfaces $y<-0.3$ and $y>0.3$, respectively. For this test case, we also want to compare the time it takes for a nonlinear scheme, namely the second-order backward differentiation formula (BDF2), and for the linear scheme introduced in our study. The BDF2 scheme is second-order accurate in time and consists in approximating the temporal derivative as in \eqref{SBDF2ThinWeakForm} and treating all other terms implicitly. This fully implicit method uses an iterative procedure to solve the nonlinear resulting system. Here, we use Newton method with a tolerance of $10^{-6}$ for the iterative solver. Numerical results for both methods are displayed in Figure~\ref{Fig3BDF2VsSBDF2}, which shows the  capability of the proposed method to simulate the macroscopic dynamics of drop splitting as with the nonlinear BDF2 scheme but with a relatively smaller computational time. In fact, we can see from Table~\ref{Tab:CPUcomparison} that the algorithm \eqref{SBDF2ThinWeakForm}-\eqref{SBDF1ThinWeakForm} is in average four times faster than the BDF2 scheme. Though the maximal time step for the proposed scheme can be smaller than for BDF methods, the required CPU time is still significantly smaller than with BDF methods.
\begin{table}[h!]
\centering
\begin{tabular}{c|c|c|c|c|c|c}
\cline{2-3}
\multicolumn{1}{c|}{}&\multicolumn{2}{|c|}{CPU (s)}\\ \cline{1-3}
\multicolumn{1}{|c|}{Final time $t_f$}&Algorithm \eqref{SBDF2ThinWeakForm}-\eqref{SBDF1ThinWeakForm}&BDF2 scheme\\ \cline{1-3}
\multicolumn{1}{|c|}{$0.018$}&$4.48$&$17.78$\\ \cline{1-3}
\multicolumn{1}{|c|}{$0.036$}&$8.36$&$34.51$\\ \cline{1-3}
\multicolumn{1}{|c|}{$0.072$}&$16.88$&$66.58$\\ \cline{1-3}
\multicolumn{1}{|c|}{$0.144$}&$33.00$&$137.66$\\ \cline{1-3}
\end{tabular}
\caption{Computation time comparison between the algorithm \eqref{SBDF2ThinWeakForm}-\eqref{SBDF1ThinWeakForm} and the BDF2 scheme for the droplet splitting test case.  $\Delta t=10^{-3}$, $\lambda=2$.}
\label{Tab:CPUcomparison}
\end{table}

\section{Conclusion}
\label{SectConclusion}
The present study provides a new fully discrete mixed finite element method to approximate fourth-order nonlinear diffusion equations. We develop a special numerical technique of approximating the nonlinear energy term that allows to design a linear scheme for our system. The performance of our scheme is tested using relevant numerical examples, and the results clearly demonstrated the ability of the proposed scheme to resolve fourth-order nonlinear diffusion equations. Manufactured and analytic solutions have been used for the numerical investigation of the accuracy of the scheme. We have performed numerical analysis and presented convergence results showing the good behavior of the numerical method with respect to mesh refinement, order of the  approximation and boundary conditions. With respect to the spatial discretization, convergence rates are obtained that match analytical error estimates available for linear problems. We also have presented a series of numerical experiments to confirm the second-order accuracy in time for both constant and non-constant mobility function. Some complex and real-world examples are presented including a model for electrowetting on dielectric and for coarsening dynamics of binary mixture.

\section*{Acknowledgement}
The first and second authors are supported by the UM6P/OCP group of Morocco. The third author is supported through a Discovery Grant of the Natural Sciences and Engineering Research Council of Canada.

\section*{References}
\bibliographystyle{plain} 
\bibliography{biblio}
\end{document}